\documentclass[a4paper,11pt]{article}
\usepackage{indentfirst,latexsym,bm,color}
\usepackage{amsmath,amssymb,amsfonts}

\usepackage[top=1in,left=1in,right=1in]{geometry}

\makeatletter

\@addtoreset{equation}{section} \makeatother
\begin{document}
\newtheorem{dingli}{Theorem}[section]
\newtheorem{dingyi}[dingli]{Definition}
\newtheorem{tuilun}[dingli]{Corollary}
\newtheorem{zhuyi}[dingli]{Remark}
\newtheorem{yinli}[dingli]{Lemma}

\title{Structure and positivity of linear maps preserving covariance under unitary evolution
\thanks{ This work was supported by NSF of
China (No: 11671242)}}
\author{Yuan Li,$^a$\thanks{E-mail address:
 liyuan0401@aliyun.com} \ \  Shuaijie Wang,$^b$\thanks{E-mail address:  wangshuaijie@xidian.edu.cn} \ \  Xiao-Ming Xu$^c$\thanks{E-mail address: xmxu@sit.edu.cn}}
\date{} \maketitle\begin{center}
\begin{minipage}{16cm}
{ \small $$ a    \ \ \    \ School \ of \ Mathematics \ and \
statistics,\ Shaanxi \ Normal \ University, $$ $$ Xi'an,\
710062,\ China. $$ }{\small $$  b \ \   \ School \ of \ Mathematics \ and \
statistics, \ Xidian \ University, $$ $$\ Xi'an, \ 710071,  \ China $$  }{ \small $$ c    \ \ \    \ School \ of \ Sciences,\ Shanghai \ Institute \ of\ Technology, $$ $$ Shanghai,\
201418,\ China. $$ }

\end{minipage}
\end{center}
 \vspace{0.05cm}
\begin{center}
\begin{minipage}{16cm}
{\small {\bf Abstract }  Let $\mathcal{H}$ be a complex finite-dimensional or infinite-dimensional separable Hilbert space, $\mathcal{B(H)}$ and $\mathcal{T(H)}$
be the Banach spaces of all bounded linear operators and of all trace class operators on $\mathcal{H},$ respectively. In this paper, we give a concrete description of the linear maps $\Phi:\mathcal{T(H)}\rightarrow \mathcal{B(H\otimes H)}$ that are continuous relative to the norm topology and covariance under unitary evolution (i.e., $\Phi(UXU^*)=(U\otimes U)\Phi(X)(U^*\otimes U^*)$ for all $X\in\mathcal{T(H)}$ and unitary operators $U\in\mathcal{B(H)}).$
Using this, we obtain the equivalent conditions for this class of maps to be self-adjoint or positive. As a corollary, we get that the virtual broadcasting map
$\mathcal{B}_{vb}:\mathcal{T(H)}\rightarrow \mathcal{B(H\otimes H)}$ with the form $\mathcal{B}_{vb}(X)=\frac{ 1}{2}[S(I\otimes X)+S(X\otimes I)]$ is uniquely determined by three conditions: covariance under unitary evolution, invariance under permutation of the copies and consistency with classical broadcasting, where $S\in\mathcal{B(H\otimes H)}$ is the swap operator. Moreover, the linear maps $\Psi:\mathcal{B(H)}\rightarrow \mathcal{B(H\otimes H)}$ that are continuous relative to the $W^*$-topology and covariance under unitary evolution are also characterized. \\}
\endabstract
\end{minipage}\vspace{0.10cm}
\begin{minipage}{16cm}
\end{minipage}
\end{center}
\begin{center} \vspace{0.01cm}
\end{center}

\section{\textbf{Introduction and Preliminary}}

Positive and completely positive linear maps had been playing important
roles in the theory of operator algebras, which reflect noncommutative order structures ([3,4,7,13,15,17,20]).
Such order structures provide basic mathematical frameworks for quantum information theory ([8-10,12,18]).
Let $\mathcal{H}$ be a complex finite-dimensional or infinite-dimensional separable Hilbert space, $\mathcal{B(H)}$
be the Banach space of all bounded linear operators on $\mathcal{H}$ with the operator norm $\|\cdot\|$ and $\mathcal{T(H)}$
be the Banach spaces of all trace class operators with  the trace norm $\|\cdot\|_1.$ As usual, $\mathcal{U(H)}$ is the group of all unitary operators on $\mathcal{H},$ $A^*$ is the adjoint operator and $|A|:=(A^*A)^\frac{1}{2}$ is the absolute value of $A.$ Also, $ran(A)$ is the range of an operator $A.$  If $A^*=A,$ then $A$ is said to be a self-adjoint operator. In particular, an operator $A\in\mathcal{B(H)}$ is called positive, if
$A\geq0,$ meaning $\langle
Ax,x\rangle \geq0 $
 for all $x\in\mathcal{H},$ where $\langle, \rangle$ is the inner product on $\mathcal{H}.$  Moreover, if $x,y\in{\mathcal{H}},$ then $xy^*$ is the one-rank operator in $\mathcal{B}(\mathcal{H})$
defined by $xy^*(z)=\langle z,y\rangle x$ for vectors $z\in\mathcal{H}.$ An operator $V\in\mathcal{B(H)}$ is called a partial isometry if $VV^*V=V.$ We always assume that the Hilbert space $\mathcal{H}$  is at least two-dimensional.

A linear map $\Phi: {\mathcal{B(H)}}\longrightarrow \mathcal{B(H\otimes H)}$ is said to be covariance under unitary evolution (unitarily equivariant) if $\Phi(UXU^*)=(U\otimes U)\Phi(X)(U^*\otimes U^*)$ for all $X\in\mathcal{B(H)}$ and $U\in\mathcal{U(H)}.$ Also, $\Phi$ is called invariance under permutation of the copies if $S\Phi(X)S=\Phi(X)$
for all $X\in\mathcal{B(H)},$ where $S\in\mathcal{B(H\otimes H)}$ is the swap operator, (i.e., $S(x\otimes y)=y\otimes x$ for all vectors $x,y\in\mathcal{H}).$

Recently, many interesting characterizations and properties of the unitarily equivariant maps  between matrix algebras were
obtained in [1,3,10,14]. In particular, Bhat in [2] proved that a linear map $\Phi$ acting on ${\mathcal{B(H)}}$ satisfies $\Phi(UXU^*)=U\Phi(X)U^*$ for all $X\in\mathcal{B(H)}$ and $U\in\mathcal{U(H)}$ if and only if there exist complex numbers $\lambda,\mu$  such that $\Phi(X)=\lambda X +\mu tr(X)I$ if $\mathcal{H}$ is finite dimensional or $\Phi(X)=\alpha X$ for a complex number $\alpha$ if $\mathcal{H}$ is infinite dimensional. More generally, $(a, b)$-unitarily equivariant linear maps between matrix algebras are study in [1]. Also, a class of linear maps that are covariant
with respect to the finite group generated by the Weyl operators (the Weyl-covariant maps) is introduced and considered in [19].

For a finite dimensional Hilbert space ${\mathcal{H}},$ a broadcasting map
 is a linear map $\Phi:\mathcal{B(H})\rightarrow \mathcal{B(H}\otimes \mathcal{H})$ satisfying
the conditions \begin{equation} tr_1[\Phi(X)]=tr_2[\Phi(X)]=X  \ \ \ \hbox{ for all }\ \ \ X\in \mathcal{B(H}),\end{equation} where $tr_{1}(A_1\otimes A_2)=tr(A_1)A_2$ and $tr_{2}(A_1\otimes A_2)=tr(A_2)A_1$ are the partial trace with respect to the first and the second Hilbert space, respectively. Note that the broadcasting condition (1.1) implies that $\Phi$ is trace preserving and the broadcasting maps have important applications in quantum information theory. However, the structural characterization of broadcasting maps seems to be difficult. We have not seen relevant results so far.

A Hermitian-preserving map that satisfies the broadcasting
condition (1.1) is referred to as a virtual broadcasting
map. It was shown in [14, Theorem 1] that the virtual broadcasting map
$\mathcal{B}_{vb}:\mathcal{B(H})\rightarrow \mathcal{B(H}\otimes \mathcal{H})$ with the form $\mathcal{B}_{vb}(X)=\frac{ 1}{2}[S(I\otimes X)+S(X\otimes I)]$ is uniquely characterized by three conditions: covariance under unitary evolution, invariance under permutation of the copies and consistency with classical broadcasting. The virtual broadcasting map $\mathcal{B}_{vb}$ has several nice properties and applications in [14, Theorem 1-4]. Moreover, the virtual broadcasting map is essential to a quantum state over time function in [6].

Let $\{e_j\}_{j=1}^s$ be
a fixed orthonormal basis for ${\mathcal{H}},$ where $s\leq\infty.$ The classical broadcasting $\mathcal{B}_{cl}: {\mathcal{B(H)}}\longrightarrow \mathcal{B(H\otimes H)}$ is defined by \begin{equation}\mathcal{B}_{cl}(X):=\sum_{i=1}^{s}\langle Xe_i,e_i\rangle e_ie_i^*\otimes e_ie_i^*,  \  \ \ X\in\mathcal{B(H)} \end{equation} We note that  $\mathcal{B}_{cl}$ is well-defined,  when ${\mathcal{H}}$ is infinite dimensional.
Moreover, a linear map $\mathcal{B}: {\mathcal{B(H)}}\longrightarrow \mathcal{B(H\otimes H)}$ is called consistency with classical broadcasting if $(\mathcal{D}\otimes\mathcal{D})\circ\mathcal{B}\circ\mathcal{D}=\mathcal{B}_{cl},$ where $\mathcal{D}$ is the diagonal operation with respect to the orthonormal basis $\{e_j\}_{j=1}^s$ (i.e.,  $\mathcal{D}(X)=\sum_{i=1}^{s}\langle Xe_i,e_i\rangle e_ie_i^*).$   Particularly, $\Phi: {\mathcal{T(H)}}\longrightarrow \mathcal{B(H\otimes H)}$ is said to be consistency with classical broadcasting if $(\mathcal{D}\otimes\mathcal{D})\circ\Phi\circ\mathcal{D}(X)=\mathcal{B}_{cl}(X)$  for all $X\in\mathcal{T(H)}.$

The aim of this paper is to consider whether there is a connection among these three conditions (covariance under unitary evolution, invariance under permutation and consistency with classical broadcasting) and how to characterize the structure of a linear map that satisfies one or two of them. Let ${\tau}_{1}$ be the topology on $\mathcal{T}(\mathcal{H})$ induced by the trace norm $\|.\|_{1}$ and ${\tau}$ be the topology on $\mathcal{B}(\mathcal{H}\otimes\mathcal{H})$ induced by the operator norm $\|.\|.$  We firstly give a concrete description of the continuous linear maps $\Phi:\mathcal{T(H)}\rightarrow \mathcal{B(H\otimes H)}$ that are covariance under unitary evolution.

{\bf Theorem 1.} Let $\Phi:(\mathcal{T}(\mathcal{H}), \tau_{1})\rightarrow (\mathcal{B}(\mathcal{H}\otimes\mathcal{H}), \tau)$ be a continuous linear  map. Then $\Phi$ is covariance under unitary evolution if and only if there exist complex numbers $\lambda_1,\lambda_2,\cdots,\lambda_6$ such that $$\Phi(X)=\lambda_1I \otimes X+\lambda_2X\otimes I+\lambda_3S(I\otimes X)+\lambda_4S(X\otimes I)+\lambda_5tr(X)I\otimes I+\lambda_6tr(X)S,$$ where $S\in\mathcal{B(H\otimes H)}$ is the swap operator.

Furthermore, the choice of complex numbers $\{\lambda_i\}_{i=1}^6$ is unique, when the dimension of $\mathcal{H}$ is at least three. However, it is not unique in the two-dimensional Hilbert space $\mathcal{H}$ (see Remak 4).

Let $\Phi:\mathcal{T(H)}\rightarrow \mathcal{B(H\otimes H)}$ be a linear map. Recall that $\Phi$ is self-adjoint (or Hermitian-preserving) if $\Phi(X^*)=\Phi(X)^*$ for all $X\in\mathcal{T(H)}.$ In particular, $\Phi$ is said to be a positive map if $\Phi(X)\geq0$ for all positive operators $X\in\mathcal{T(H)}.$
Based on the above Theorem 1, we derive some equivalent conditions under which $\Phi$ becomes a self-adjoint map or a positive map.

{\bf Theorem 2.} Let $S\in\mathcal{B(H\otimes H)}$ be the swap operator and $\Phi:(\mathcal{T}(\mathcal{H}), \tau_{1})\rightarrow (\mathcal{B}(\mathcal{H}\otimes\mathcal{H}), \tau)$ be a continuous linear  map, which  is covariance under unitary evolution.  

$(a)$ $\Phi$ is a self-adjoint map if and only if there exist real numbers $\mu_i,$ $i=1,2,5,6$ and complex numbers $\mu_3=\overline{\mu_4}$ such that \begin{equation}\Phi(X)=\mu_1I \otimes X+\mu_2X\otimes I+\mu_3S(I\otimes X)+\mu_4S(X\otimes I)+\mu_5tr(X)I\otimes I+\mu_6tr(X)S.\end{equation}

$(b)$ $\Phi$ is a positive map if and only if $\Phi$ has the form (1.3) and the coefficients satisfy  $\sum_{i=1}^{6}\mu_i\geq0,$
 $\left(\begin{array}{cc}
\mu_1+\mu_5&\overline{\mu_3}+\mu_6\\ \mu_3+\mu_6&\mu_2+\mu_5\end{array}\right)\geq0$
and $$
    \left\{
    \begin{array}{l}
     \mu_5\geq|\mu_6|,      \ \ \ \ \ \ \ \ \ \ \ if \  \ dim \mathcal{H}\geq 3  \ \  \ \ \\
     \mu_5+\mu_6\geq0,    \ \  \ \ \ \ \ if \ \ dim \mathcal{H}=2.\ \ \ \ \ \   \ \\
        \end{array}
 \right.$$

As a corollary, we obtain in particular that the virtual broadcasting map
$\mathcal{B}_{vb}:\mathcal{T(H)}\rightarrow \mathcal{B(H\otimes H)}$ with the form $\mathcal{B}_{vb}(X)=\frac{ 1}{2}[S(I\otimes X)+S(X\otimes I)]$ is uniquely characterized by the above three conditions: covariance under unitary evolution, invariance under permutation of the copies and consistency with classical broadcasting. Moreover, we consider the
specific structure of the linear maps $\Psi:\mathcal{B(H)}\rightarrow \mathcal{B(H\otimes H)}$ that are continuous in the $W^*$-topology and covariance under unitary evolution. In this case, we get that $\Psi$ is a completely positive map if and only if it is a positive map (see Corollary 11).

{\bf Theorem 3.} Let ${\mathcal{H}}$ be an infinite-dimensional separable Hilbert space and $\Psi:\mathcal{B(H)}\rightarrow \mathcal{B(H\otimes H)}$ be a continuous linear map in the $W^*$-topology. Then the following statements are equivalent:

$(a)$ $\Psi$ is covariance under unitary evolution;

$(b)$ $\Psi(Uee^*U^*)=(U\otimes U)\Psi(ee^*)(U^*\otimes U^*)$ for all unit vectors $e\in\mathcal{H}$ and $U\in\mathcal{U(H)};$

$(c)$ There exist complex numbers $\lambda_i$ for $i=1,2,3,4$ such that $$\Psi(X)=\lambda_1I \otimes X+\lambda_2X\otimes I+\lambda_3S(I\otimes X)+\lambda_4S(X\otimes I),$$ where $S\in\mathcal{B(H\otimes H)}$ is the swap operator.

Let $m$ be a positive integer. To get a concrete description of the continuous linear maps $\Phi:\mathcal{T(H)}\rightarrow \mathcal{B(H}^{\otimes m})$ that are covariance under unitary evolution $(\Phi(UXU^*)=U^{\otimes m}\Phi(X)(U^*)^{\otimes m}$ for all $X\in\mathcal{T(H)}$ and $U\in\mathcal{U(H)}),$ we need to use the Schur-Weyl duality.

Let $\mathcal{C}^d$ denote a $d$-dimensional complex Hilbert space. Consider the following representation of the unitary group $\mathcal{U(C}^d)$ on $(\mathcal{C}^d)^{\otimes m}.$ For any $U\in\mathcal{U(C}^d),$ we define $\Delta(U)\in \mathcal{B((C}^d)^{\otimes m})$ by
$$\Delta(U)(x_1\otimes x_2\otimes\cdots\otimes x_m)=Ux_1\otimes Ux_2\otimes\cdots\otimes U x_m.$$
Obviously, $\Delta(U)=U\otimes U\otimes\cdots\otimes U=U^{\otimes m}.$ Let $\mathbb{S}_m$
be the symmetric group of degree $m.$ Consider the canonical representation of the symmetric group $\mathbb{S}_m$ on $(\mathcal{C}^d)^{\otimes m}.$
That is, $$\Gamma(s)(x_1\otimes x_2\otimes\cdots\otimes x_m)=x_{s^{-1}(1)}\otimes x_{s^{-1}(2)}\otimes\cdots\otimes x_{s^{-1}(m)},$$ for $s\in\mathbb{S}_m.$ It is clear that $\Gamma(s)$ are all unitary operators. We denote the groups $$\Delta(\mathcal{U(C}^d))=\{U^{\otimes m}:\ U \in \mathcal{U(C}^d)\} \ \hbox{ and }\ \Gamma(\mathbb{S}_m)=\{\Gamma(s): \ s \in \mathbb{S}_m\}.$$ Then Schur-Weyl duality describes the commutativity of the algebra generated by the two sets above as follows.

{\bf Theorem (Schur-Weyl duality).} ([11]) The following two algebras are commutants of one another in  $\mathcal{B((C}^d)^{\otimes m}).$

$(a)$ $Alg\{\Delta(\mathcal{U(C}^d))\},$ the complex algebra spanned by $\Delta(\mathcal{U(C}^d)).$

$(b)$ $Alg\{\Gamma(\mathbb{S}_m)\},$ the complex algebra spanned by $\Gamma(\mathbb{S}_m).$

In the following, we give a concrete description of a continuous linear map $\Phi:\mathcal{T(H)}\rightarrow \mathcal{B(H}^{\otimes m})$ that is covariance under unitary evolution.

{\bf Theorem 4.} Let $dim\mathcal{H}\geq3,$ and let  $\Phi:(\mathcal{T(H)},\|\cdot\|_1)\rightarrow (\mathcal{B(H}^{\otimes m}),\|\cdot\|)$ be a continuous linear map. Then
$\Phi$ is covariance under unitary evolution if and only if there exist complex numbers $\lambda_{ij}$ for $i=1,2,\cdots,m!$ and $j=1,2,\cdots,m+1$ such that $$\Phi(X)=\sum\limits_{i=1}^{m!}\sum\limits_{j=1}^{m+1}\lambda_{ij}\Gamma(s_{i})\Phi_j(X),$$ where $s_k\neq s_i\in\mathbb{S}_m$ for $k\neq i=1,2,\cdots,m!,$ $\Phi_{1}(X) =tr(X)I^{\otimes m}$ and
 $\Phi_k(X)=I^{\otimes (k-2)}\otimes X\otimes I^{\otimes (m+1-k)}$ for $2\leq k\leq m+1.$
 Moreover, if $dim\mathcal{H}\geq m+1,$ then the choice of above complex numbers $\lambda_{ij}$   is unique.

Furthermore, we consider how to characterize the completely bounded norm $\|\Psi\|_{cb}$ of $\Psi,$ where $\Psi:\mathcal{B(H)}\rightarrow \mathcal{B(H\otimes H)}$ is a continuous linear map in the $W^*$-topology and covariance under unitary evolution.

\section{Corollaries and Proofs of Theorems 1 and 2}

The following result is well-known by using Schur-Weyl duality ([11,21]), which is a powerful tool in representation theory and has
many applications to quantum information.

{\bf Lemma 1.} Let $\mathcal{H}$ be a finite dimensional Hilbert space and $T\in\mathcal{B(H\otimes H)}.$ Then $T(U\otimes U)=(U\otimes U)T$ for all $U\in\mathcal{U(H)}$ if and only if  there exist complex numbers $\lambda$ and $\mu$ such that $T=\lambda I\otimes I+\mu S,$ where $S\in\mathcal{B(H\otimes H)}$ is the swap operator.

In the following, we shall extend Lemma 1 to a separable Hilbert space $\mathcal{H}.$ Recall that a sequence $\{A_n\}_{n=1}^\infty\subseteq\mathcal{B(H)}$ converges to $A$ in the weak*-topology,
denoted $W^*$-$\lim\limits_{n\rightarrow\infty}A_n=A,$ if $\lim\limits_{n\rightarrow\infty}tr(A_nT)=tr(AT)$ for all $T\in\mathcal{T(H)}.$ In particular,  if $\{e_j\}_{j=1}^\infty$ is an orthonormal basis for ${\mathcal{H}}$ and
$P_n=\sum_{i=1}^{n}  e_ie_i^*,$ then $W^*$-$\lim\limits_{n\rightarrow\infty}P_nAP_n=A$  for any $A\in\mathcal{B(H)}$ (see [5, Chapters 1.6 and 2.2]).

{\bf Lemma 2.} Let ${\mathcal{H}}$ be an infinite-dimensional separable Hilbert space and  $T\in\mathcal{B(H\otimes H)}.$ Then $T(U\otimes U)=(U\otimes U)T$ for all $U\in\mathcal{U(H)}$ if and only if  there exist complex numbers $\lambda$ and $\mu$ such that $T=\lambda I\otimes I+\mu S,$ where $S\in\mathcal{B(H\otimes H)}$ is the swap operator.

{\bf Proof.} Sufficiency is clear.

Necessity. Let $\{e_j\}_{j=1}^\infty$ be
an orthonormal basis for ${\mathcal{H}},$ and let $P_n=\sum\limits_{i=1}^{n} e_ie_i^*$ and $$T_n=(P_n\otimes P_n)T(P_n\otimes P_n)\mid_{ran(P_n)\otimes ran(P_n)}.$$ Then $W^*$-$\lim\limits_{n\rightarrow\infty}P_n =I,$ and since $T(U\otimes U)=(U\otimes U)T$ for all $U\in\mathcal{U(H)},$ we have $$T_n(V\otimes V)=(V\otimes V)T_n  \ \ \hbox { for all }\ \ V\in\mathcal{U}(ran(P_n)).$$ Now applying Lemma 1 to $T_{n}$ for $n=2, 3, ...,$ we can find complex numbers $\lambda_n$ and $\mu_n$ such that $$T_n=\lambda_n P_n\otimes P_n\mid_{ran(P_n)\otimes ran(P_n)}+\mu_n S_n,$$ where $S_n\in\mathcal{B}(ran(P_n)\otimes ran(P_n))$ is the swap operator.
 
We claim that $\lambda_2=\lambda_3=\cdots$ and $\mu_2=\mu_3=\cdots.$
If $2\leq m<n,$ we have $$\begin{array}{rcl}\lambda_m P_m\otimes P_m+\mu_m S(P_m\otimes P_m)&=&(P_m\otimes P_m)T(P_m\otimes P_m)\\&=&(P_m\otimes P_m)(P_n\otimes P_n)T(P_n\otimes P_n)(P_m\otimes P_m).\end{array}$$ Since the fact that $$S_n=S(P_n\otimes P_n)\mid_{ran(P_n)\otimes ran(P_n)}$$ gives
$$(P_n\otimes P_n)T(P_n\otimes P_n)=\lambda_n P_n\otimes P_n+\mu_n S(P_n\otimes P_n),$$ which implies $$\begin{array}{rcl}(P_m\otimes P_m)[\lambda_n P_n\otimes P_n+\mu_n S(P_n\otimes P_n)](P_m\otimes P_m)=\lambda_n P_m\otimes P_m+\mu_n S(P_m\otimes P_m).\end{array}$$ It entails that $\lambda_m=\lambda_n$ and $\mu_m=\mu_n.$ 
Therefore, $$\begin{array}{rcl}T&=&W^*\hbox{-}\lim\limits_{n\rightarrow\infty}(P_n\otimes P_n)T(P_n\otimes P_n)\\&=&W^*\hbox{-}\lim\limits_{n\rightarrow\infty}\lambda_n P_n\otimes P_n+W^*\hbox{-}\lim\limits_{n\rightarrow\infty}\mu_n S(P_n\otimes P_n)\\&=&\lambda I\otimes I+\mu S.\end{array}$$ $\Box$

{\bf Proof of Theorem 1.} Sufficiency is easy to verify.

Necessity.
Let $e\in \mathcal{H}$ be a fixed unit vector. The first step is to establish the following.

{\bf Claim 1.} There exist complex numbers $\lambda_i$ for $i=0,1,\cdots,6$ such that $$\Phi(ee^*)=\lambda_0ee^*\otimes ee^*+\lambda_1I \otimes ee^*+\lambda_2ee^*\otimes I+\lambda_3S(I\otimes ee^*)+\lambda_4S( ee^*\otimes I)+\lambda_5I\otimes I+\lambda_6S.$$
  Let 
$$\mathcal{U}_e:=\{ee^*+ V: \ V \hbox{ is a partial isometry with } ran(V)=ran(V^*)=ran(I-ee^*) \}.$$
For any $U\in \mathcal{U}_e,$ it is easy to see that $U$ is a unitary operator with $Uee^*=ee^*U.$ Then $$\Phi(ee^*)=\Phi(Uee^*U^*)=(U\otimes U)\Phi(ee^*)(U^*\otimes U^*),$$ and multiplying by $U{\otimes}U$ gives $(U\otimes U)\Phi(ee^*)=\Phi(ee^*)(U\otimes U).$ Notice also that 
 $$U\otimes U=(ee^*+ V)\otimes(ee^*+ V)=ee^*\otimes ee^*+ ee^*\otimes V+ V\otimes ee^*+V\otimes V$$ for some partial isometry $V$ with $ran(V)=ran(V^*)=ran(I-ee^*),$  it becomes  \begin{equation}(ee^*\otimes ee^*+ ee^*\otimes V+ V\otimes ee^*+V\otimes V)\Phi(ee^*)=\Phi(ee^*)( ee^*\otimes ee^*+ ee^*\otimes V+ V\otimes ee^*+V\otimes V).\end{equation}  Thus \begin{equation}(ee^*\otimes ee^*+ ee^*\otimes V+ V\otimes ee^*+V\otimes V)\Phi(ee^*)(ee^*\otimes ee^*)=\Phi(ee^*)(ee^*\otimes ee^*),\end{equation} and so
 \begin{equation}(ee^*\otimes V+ V\otimes ee^*+V\otimes V)\Phi(ee^*)(ee^*\otimes ee^*)=\Phi(ee^*)(ee^*\otimes ee^*)-(ee^*\otimes ee^*)\Phi(ee^*)(ee^*\otimes ee^*).\end{equation}
Replacing $V$ by $-V$ in equation (2.3),  we also know that  \begin{equation}(-ee^*\otimes V-V\otimes ee^*+V\otimes V)\Phi(ee^*)(ee^*\otimes ee^*)=\Phi(ee^*)(ee^*\otimes ee^*)-(ee^*\otimes ee^*)\Phi(ee^*)(ee^*\otimes ee^*).\end{equation}
Combining equations (2.3) and (2.4), a little arithmetic implies   \begin{equation}(V\otimes V)\Phi(ee^*)(ee^*\otimes ee^*)=\Phi(ee^*)(ee^*\otimes ee^*)-(ee^*\otimes ee^*)\Phi(ee^*)(ee^*\otimes ee^*) \end{equation} and \begin{equation}(ee^*\otimes V+V\otimes ee^*)\Phi(ee^*)(ee^*\otimes ee^*)=0. \end{equation}
But, using $\sqrt{-1}V$ instead of $V$ in equation (2.5) gives \begin{equation}-(V\otimes V)\Phi(ee^*)(ee^*\otimes ee^*)=\Phi(ee^*)(ee^*\otimes ee^*)-(ee^*\otimes ee^*)\Phi(ee^*)(ee^*\otimes ee^*), \end{equation} so \begin{equation}(V\otimes V)\Phi(ee^*)(ee^*\otimes ee^*)=0=\Phi(ee^*)(ee^*\otimes ee^*)-(ee^*\otimes ee^*)\Phi(ee^*)(ee^*\otimes ee^*).\end{equation}

In a similar way, it also follows from (2.1) that \begin{equation}(ee^*\otimes ee^*)\Phi(ee^*)(ee^*\otimes ee^*+ ee^*\otimes V+ V\otimes ee^*+V\otimes V)=(ee^*\otimes ee^*)\Phi(ee^*),\end{equation}
which implies \begin{equation}(ee^*\otimes ee^*)\Phi(ee^*)(V\otimes V)=0=(ee^*\otimes ee^*)\Phi(ee^*)-(ee^*\otimes ee^*)\Phi(ee^*)(ee^*\otimes ee^*).\end{equation}
Thus \begin{equation}\Phi(ee^*)(ee^*\otimes ee^*)=(ee^*\otimes ee^*)\Phi(ee^*)=(ee^*\otimes ee^*)\Phi(ee^*)(ee^*\otimes ee^*),\end{equation} which yields that there exists a complex number  $\mu_1$ such that
\begin{equation}\Phi(ee^*)(ee^*\otimes ee^*)=(ee^*\otimes ee^*)\Phi(ee^*)=\mu_1ee^*\otimes ee^*.\end{equation}
Moreover, equations (2.1) and (2.12) imply \begin{equation}(ee^*\otimes V+ V\otimes ee^*+V\otimes V)\Phi(ee^*)=\Phi(ee^*)(ee^*\otimes V+ V\otimes ee^*+V\otimes V).\end{equation}
Replacing $V$ by $-V$ in equation (2.13),  we have that \begin{equation}(-ee^*\otimes V- V\otimes ee^*+V\otimes V)\Phi(ee^*)=\Phi(ee^*)(-ee^*\otimes V- V\otimes ee^*+V\otimes V),\end{equation}
and so \begin{equation}(ee^*\otimes V+ V\otimes ee^*)\Phi(ee^*)=\Phi(ee^*)(ee^*\otimes V+ V\otimes ee^*) \end{equation} and \begin{equation}(V\otimes V)\Phi(ee^*)=\Phi(ee^*)(V\otimes V).\end{equation}

Denote by $P:=I-ee^*.$ It is clear that $PV=VP=V.$ By equation (2.15), we see that  
\begin{equation}( V \otimes ee^*)\Phi(ee^*)( P \otimes ee^*)=( P\otimes ee^*)\Phi(ee^*)(V\otimes ee^*)\end{equation} and \begin{equation}(ee^*\otimes V)\Phi(ee^*)(ee^*\otimes P)=(ee^*\otimes P)\Phi(ee^*)(ee^*\otimes V).\end{equation} Thus there exist two complex number  $\mu_2$ and $\mu_3$ such that \begin{equation}(P\otimes ee^*)\Phi(ee^*)(P \otimes ee^*)= \mu_2P \otimes ee^*\end{equation} and \begin{equation}(ee^*\otimes P )\Phi(ee^*)(ee^*\otimes P )= \mu_3 ee^*\otimes P, \end{equation} since $V|_{ran(I-ee^*)}$ is an arbitrary unitary operator acting on the subspace of $ran(I-ee^*).$ Furthermore, we conclude from equation (2.15) that \begin{equation}
(ee^*\otimes V)\Phi(ee^*)(P\otimes ee^*)=(ee^*\otimes P)\Phi(ee^*)( V\otimes ee^*)\end{equation}
and \begin{equation}( V \otimes ee^*)\Phi(ee^*)(ee^*\otimes P)=( P\otimes ee^*)\Phi(ee^*)( ee^*\otimes V).\end{equation}  Since $S\in\mathcal{B(H\otimes H)}$ is the swap operator,  
equation (2.21) implies that \begin{equation}\begin{array}{rcl}
(V\otimes ee^*)S\Phi(ee^*)(P\otimes ee^*)&=&S(ee^*\otimes V)\Phi(ee^*)(P\otimes ee^*)\\&=&S(ee^*\otimes P)\Phi(ee^*)( V\otimes ee^*)\\&=&(P\otimes ee^*)S\Phi(ee^*)( V\otimes ee^*).\end{array}\end{equation}
So there exists a complex number  $\mu_4$ such that \begin{equation}(P\otimes ee^*)S\Phi(ee^*)(P \otimes ee^*)= \mu_4P \otimes ee^*,\end{equation} and thus \begin{equation}(ee^*\otimes P)\Phi(ee^*)(P \otimes ee^*)=\mu_4S(P\otimes ee^*).\end{equation}
 In a similar way, it follows from equation (2.22) that there exists a complex number $\mu_5$ such that
\begin{equation}(P \otimes ee^*)\Phi(ee^*)(ee^*\otimes P)=\mu_5S(ee^*\otimes P).\end{equation}

Clearly, equation (2.16) also implies \begin{equation}(V\otimes V)[(P\otimes P)\Phi(ee^*)(P\otimes P)]=[(P\otimes P)\Phi(ee^*)(P\otimes P)](V\otimes V).\end{equation} Then we conclude from Lemma 2 (or Lemma 1) that there exist two complex numbers $\mu_6$ and $\mu_7$ such that \begin{equation}(P\otimes P)\Phi(ee^*)(P\otimes P)=\mu_6P\otimes P+\mu_7S(P\otimes P).\end{equation} It is obvious that \begin{equation}(ee^*\otimes P+ P\otimes ee^*)\Phi(ee^*)=\Phi(ee^*)(ee^*\otimes P+ P\otimes ee^*)\end{equation} and  \begin{equation}(P\otimes P)\Phi(ee^*)=\Phi(ee^*)(P\otimes P)\end{equation}
follow from equations (2.15) and (2.16). Moreover, it is easy to check that $$S(ee^*\otimes ee^*)(x\otimes y)=S[(e^*x) e\otimes(e^*y)e]=(e^*x)(e^*y) e\otimes e=(ee^*\otimes ee^*)(x\otimes y)$$ for all $x,y\in\mathcal{H},$ so \begin{equation}S(ee^*\otimes ee^*)=ee^*\otimes ee^*.\end{equation}
Combining equations (2.12),(2.19)-(2.20), (2.25)-(2.26) and (2.28)-(2.31), we get that \begin{equation}\begin{array}{rcl}
\Phi(ee^*)&=&\Phi(ee^*)(ee^*\otimes ee^*)+\Phi(ee^*)(ee^*\otimes P+P\otimes ee^*)+(P\otimes P)\Phi(ee^*)(P\otimes P)\\&=&\mu_1ee^*\otimes ee^*+\mu_2P \otimes ee^*+\mu_3 ee^*\otimes P+\mu_4S(P\otimes ee^*)\\&&+\mu_5S( ee^*\otimes P)+\mu_6P\otimes P+\mu_7S(P\otimes P)\\&=&(\mu_1-\mu_6-\mu_7)ee^*\otimes ee^*+(\mu_2-\mu_6)P \otimes ee^*+(\mu_3-\mu_6) ee^*\otimes P\\&&+(\mu_4-\mu_7)S(P\otimes ee^*)+(\mu_5-\mu_7)S( ee^*\otimes P)+\mu_6I\otimes I+\mu_7S\\&=&(\mu_1-\sum_{j=2}^{5}\mu_j+\mu_6+\mu_7)ee^*\otimes ee^*+(\mu_2-\mu_6)I \otimes ee^*+(\mu_3-\mu_6) \\&&ee^*\otimes I+(\mu_4-\mu_7)S(I\otimes ee^*)+(\mu_5-\mu_7)S( ee^*\otimes I)+\mu_6I\otimes I+\mu_7S
\\&=&\lambda_0ee^*\otimes ee^*+\lambda_1I \otimes ee^*+\lambda_2ee^*\otimes I+\lambda_3S(I\otimes ee^*)+\lambda_4S( ee^*\otimes I)\\&&+\lambda_5I\otimes I+\lambda_6S,\end{array}\end{equation}   where $$\lambda_0:=\mu_1-\sum_{j=2}^{5}\mu_j+\mu_6+\mu_7, \ \
\lambda_j:=\mu_{j+1}-\mu_6,\ \ \ \lambda_{j+2}:=\mu_{j+3}-\mu_7 \ \hbox{ and } \lambda_{j+4}=\mu_{j+5}$$ for $j=1,2.$ This establishes the above claim 1.

{\bf Claim 2.}  Let $f$ be any unit vector in $\mathcal{H}.$ Then $$\Phi(ff^*)=\lambda_1I \otimes ff^*+\lambda_2ff^*\otimes I+\lambda_3S(I\otimes ff^*)+\lambda_4S( ff^*\otimes I)+\lambda_5I\otimes I+\lambda_6S.$$

Suppose that $U_{f}\in\mathcal{B}(\mathcal{H})$ is a unitary operator  with $U_fe=f.$ It is obvious that \begin{equation}\begin{array}{rcl}\Phi(ff^*)&=&\Phi(U_fee^*U_f^*)\\&=&(U_f\otimes U_f)\Phi(ee^*)(U_f\otimes U_f)^*\\&=&(U_f\otimes U_f)[\lambda_0ee^*\otimes ee^*+\lambda_1I \otimes ee^*+\lambda_2ee^*\otimes I+\lambda_3S(I\otimes ee^*)\\&&+\lambda_4S( ee^*\otimes I)+\lambda_5I\otimes I+\lambda_6S](U_f\otimes U_f)^*\\&=&  \lambda_0ff^*\otimes ff^*+\lambda_1I \otimes ff^*+\lambda_2ff^*\otimes I+\lambda_3S(I\otimes ff^*)\\&&+\lambda_4S(ff^*\otimes I)+\lambda_5I\otimes I+\lambda_6S.\end{array}\end{equation}

Let $ \widetilde{e}$ be a unit vector with $\widetilde{e}\bot e.$ Setting $g=\frac{1}{\sqrt{2}}e+\frac{1}{\sqrt{2}}\widetilde{e}$ and $h=\frac{1}{\sqrt{2}}e-\frac{1}{\sqrt{2}}\widetilde{e},$ we get that $$gg^*+hh^*=ee^*+\widetilde{e}\widetilde{e}^*$$ and 
\begin{equation}\Phi(gg^*)+\Phi(hh^*)=\Phi(ee^*)+\Phi(\widetilde{e}\widetilde{e}^*).\end{equation} It follows from equations (2.33) and (2.34) that \begin{equation}\begin{array}{rcl}
\lambda_0gg^*\otimes gg^*+\lambda_0hh^*\otimes hh^*=\lambda_0ee^*\otimes ee^*+\lambda_0\widetilde{e}\widetilde{e}^*\otimes \widetilde{e}\widetilde{e}^*,\end{array}\end{equation}   which yields \begin{equation}\begin{array}{rcl}
0&=&\lambda_0(gg^*\otimes gg^*+hh^*\otimes hh^*)(g\otimes h)\\&=&\lambda_0(ee^*\otimes ee^*+\widetilde{e}\widetilde{e}^*\otimes \widetilde{e}\widetilde{e}^*)(g\otimes h)=\frac{\lambda_0}{2}(e\otimes e-\widetilde{e}\otimes\widetilde{e}),\end{array}\end{equation}
and hence $\lambda_0=0.$ Thus  equation (2.33) yields that \begin{equation}\Phi(ff^*)=\lambda_1I \otimes ff^*+\lambda_2ff^*\otimes I+\lambda_3S(I\otimes ff^*)+\lambda_4S(ff^*\otimes I)+\lambda_5I\otimes I+\lambda_6S,\end{equation} and the claim 2 is proved.

For any positive operator $A\in\mathcal{T(H)},$ we get from the spectral decomposition theorem that $A=\sum_{j=1}^{\infty}t_jf_jf_j^*,$ where $\{f_j\}_{j=1}^\infty$ is
an orthonormal basis of ${\mathcal{H}}$ and $t_j\geq t_{j+1}\geq0$ are eigenvalues of $A.$ Denoting by $A_n=\sum_{j=1}^{n}t_jf_jf_j^*,$ we have  $$\lim\limits_{n\rightarrow\infty}\|A_n-A\|_1=\lim\limits_{n\rightarrow\infty}\sum_{j=n+1}^{\infty}t_j=0.$$   Then
$$\begin{array}{rl}&\Phi(A)\\=&\lim\limits_{n\rightarrow\infty}\Phi(A_n)\\=&\lim\limits_{n\rightarrow\infty}\sum_{j=1}^{n}
t_j\Phi(f_jf_j^*)\\=&\lim\limits_{n\rightarrow\infty}\sum_{j=1}^{n}t_j[\lambda_1I \otimes f_jf_j^*+\lambda_2f_jf_j^*\otimes I+\lambda_3S(I\otimes f_jf_j^*)+\lambda_4S(f_jf_j^*\otimes I)+\lambda_5I\otimes I+\lambda_6S]\\=&\lim\limits_{n\rightarrow\infty}[\lambda_1I \otimes A_n+\lambda_2A_n\otimes I+\lambda_3S(I\otimes A_n)+\lambda_4S(A_n\otimes I)+\lambda_5tr(A_n) I\otimes I+\lambda_6tr(A_n)S]\\=&\lambda_1I \otimes A+\lambda_2A\otimes I+\lambda_3S(I\otimes A)+\lambda_4S(A\otimes I)+\lambda_5tr(A) I\otimes I+\lambda_6tr(A)S\end{array}$$ follows from the facts $$\lim\limits_{n\rightarrow\infty}\|A_n-A\|\leq\lim\limits_{n\rightarrow\infty}\|A_n-A\|_1=0$$  and $$ \lim\limits_{n\rightarrow\infty}|tr(A_n)- tr(A)|=\lim\limits_{n\rightarrow\infty}\sum_{j=n+1}^{\infty}\lambda_j=0.$$  If $A\in\mathcal{T(H)}$ is a self-adjoint operator, then $A=A^+-A^-$ where $A^+=\frac{|A|+A}{2}$ and $A^-=\frac{|A|-A}{2}$ are positive. Thus $$\begin{array}{rcl}&&\Phi(A)\\&=&\Phi(A^+)-\Phi(A^-)\\&=&\lambda_1I \otimes A^++\lambda_2A^+\otimes I+\lambda_3S(I\otimes A^+)+\lambda_4S(A^+\otimes I)+\lambda_5tr(A^+) I\otimes I+\lambda_6tr(A^+)S
\\&&-\lambda_1I \otimes A^--\lambda_2A^-\otimes I-\lambda_3S(I\otimes A^-)-\lambda_4S(A^-\otimes I)-\lambda_5tr(A^-) I\otimes I-\lambda_6tr(A^-)S\\&=&\lambda_1I \otimes A+\lambda_2A\otimes I+\lambda_3S(I\otimes A)+\lambda_4S(A\otimes I)+\lambda_5tr(A) I\otimes I+\lambda_6tr(A)S.\end{array}$$
For a general $X\in\mathcal{T}\mathcal{(H)},$ we can write
$X=Re(X)+\sqrt{-1}Im(X),$ where $Re(X)=\frac{X+X^*}{2}$ and $Im(X)=\frac{X-X^*}{2\sqrt{-1}}$ are self-adjoint operators. Then $$\begin{array}{rcl}\Phi(X)&=&\Phi(Re(X))+\sqrt{-1}\Phi(Im(X))\\&=&\lambda_1I \otimes X+\lambda_2X\otimes I+\lambda_3S(I\otimes X)+\lambda_4S(X\otimes I)+\lambda_5tr(X) I\otimes I+\lambda_6tr(X)S.\end{array}$$
$\Box$

{\bf Proof of Theorem 2.} $(a)$ Sufficiency is easy to verify.

Necessity.
Using Theorem 1, we get that there exist complex numbers $\lambda_i$ for $i=1,2,\cdots,6$ such that \begin{equation}\Phi(X)=\lambda_1I \otimes X+\lambda_2X\otimes I+\lambda_3S(I\otimes X)+\lambda_4S(X\otimes I)+\lambda_5tr(X)I\otimes I+\lambda_6tr(X)S.\end{equation} Here $S$ is the swap operator. Let $e_1$ and $e_2$ be orthogonal unit vectors.  
Then $$\begin{array}{rcl}[\Phi(e_1e_2^*)]^*&=&[\lambda_1I \otimes e_1e_2^*+\lambda_2e_1e_2^*\otimes I+\lambda_3S(I\otimes e_1e_2^*)+\lambda_4S(e_1e_2^*\otimes I)]^*\\&=&\overline{\lambda_1}I \otimes e_2e_1^*+\overline{\lambda_2}e_2e_1^*\otimes I+\overline{\lambda_3}(I\otimes e_2 e_1^*)S+\overline{\lambda_4}(e_2 e_1^*\otimes I)S\\&=&\overline{\lambda_1}I \otimes e_2e_1^*+\overline{\lambda_2}e_2e_1^*\otimes I+\overline{\lambda_3}S(e_2 e_1^*\otimes I)+\overline{\lambda_4}S(I\otimes e_2 e_1^*)\end{array}$$ and $$\Phi(e_2e_1^*)=\lambda_1I \otimes e_2 e_1^*+\lambda_2e_2e_1^*\otimes I+\lambda_3S(I\otimes e_2e_1^*)+\lambda_4S(e_2e_1^*\otimes I).$$
 Notice also that $(\Phi(e_{1}e_{2}^*))^*=\Phi(e_{2}e_{1}^*)$ because $\Phi$ is self-adjoint, we see that  \begin{equation}(\lambda_1-\overline{\lambda_1})I \otimes e_2e_1^*+(\lambda_2-\overline{\lambda_2})e_2e_1^*\otimes I+(\lambda_3-\overline{\lambda_4})S(I\otimes e_2 e_1^*)+(\lambda_4-\overline{\lambda_3})S(e_2 e_1^*\otimes I)=0.\end{equation} It entails that \begin{equation}\begin{array}{rl}&(\lambda_2-\overline{\lambda_2})e_2e_2^*\otimes e_2e_1^*+(\lambda_4-\overline{\lambda_3})(e_2e_2^*\otimes e_2e_1^*)\\=&(\lambda_2-\overline{\lambda_2})e_2e_2^*\otimes e_2e_1^*+(\lambda_4-\overline{\lambda_3})S(e_2e_2^*\otimes e_2e_1^*)\\=&[(\lambda_1-\overline{\lambda_1})I \otimes e_2e_1^*+(\lambda_2-\overline{\lambda_2})e_2e_1^*\otimes I+(\lambda_3-\overline{\lambda_4})S(I\otimes e_2 e_1^*)\\&+(\lambda_4-\overline{\lambda_3})S(e_2 e_1^*\otimes I)](e_1 e_2^*\otimes e_2e_1^*)\\=&0,\end{array}\end{equation}
 and so \begin{equation}\lambda_2-\overline{\lambda_2}+\lambda_4-\overline{\lambda_3}=0.\end{equation}
In a similar way, we obtain that \begin{equation}\begin{array}{rl}&(\lambda_1-\overline{\lambda_1})e_2e_1^*\otimes e_2e_2^*+(\lambda_3-\overline{\lambda_4})(e_2e_1^*\otimes e_2e_2^*)\\=&[(\lambda_1-\overline{\lambda_1})I \otimes e_2e_1^*+(\lambda_2-\overline{\lambda_2})e_2e_1^*\otimes I+(\lambda_3-\overline{\lambda_4})S(I\otimes e_2 e_1^*)\\&+(\lambda_4-\overline{\lambda_3})S(e_2 e_1^*\otimes I)](e_2e_1^* \otimes e_1 e_2^*)\\=&0,\end{array}\end{equation}
and hence, \begin{equation}\lambda_1-\overline{\lambda_1}+\lambda_3-\overline{\lambda_4}=0.\end{equation}

 {\bf Case 1.} If $dim \mathcal{H}\geq 3,$ then we can find a unit vector $e_{3}$ in $\mathcal{H}$ such that  $e_3\bot e_2$ and $e_3\bot e_1.$ It follows from equation (2.39) that \begin{equation}\begin{array}{rl}&(\lambda_1-\overline{\lambda_1})e_3\otimes e_2+(\lambda_4-\overline{\lambda_3})e_2\otimes e_3\\=&[(\lambda_1-\overline{\lambda_1})I \otimes e_2e_1^*+(\lambda_2-\overline{\lambda_2})e_2e_1^*\otimes I+(\lambda_3-\overline{\lambda_4})S(I\otimes e_2 e_1^*)\\&+(\lambda_4-\overline{\lambda_3})S(e_2 e_1^*\otimes I)](e_3\otimes e_1)\\=&0.\end{array}\end{equation}
It yields $\lambda_1=\overline{\lambda_1}$ and $\lambda_4=\overline{\lambda_3},$ and by equation (2.41), we have  $\lambda_2=\overline{\lambda_2}.$ 
 Since $\Phi$ is a self-adjoint map, equation (2.38) induces that $$\Lambda(X)=\lambda_5tr(X)I\otimes I+\lambda_6tr(X)S$$ is also a self-adjoint map. So  $\lambda_5=\overline{\lambda_5}$ and $\lambda_6=\overline{\lambda_6}.$ Setting $\mu_i=\lambda_i$ for $i=1,2,\cdots,6,$ we get the desired result.

 {\bf Case 2.} If $dim \mathcal{H}=2,$ then $\mathcal{H}=\bigvee\{e_1,e_2\}.$
It follows from equations (2.41) and (2.43) that \begin{equation}\lambda_2-\overline{\lambda_2}=\overline{\lambda_3}-\lambda_4=\lambda_1-\overline{\lambda_1}.\end{equation}
Denote by $\lambda_i=\alpha_i+\beta_i\sqrt{-1}$ for $i=1,2,\cdots,6.$ Then equation (2.45) implies that \begin{equation}2\beta_1=2\beta_2=-(\beta_3+\beta_4)\ \ \hbox{ and } \  \ \alpha_3=\alpha_4.\end{equation}
Clearly, $\Phi(e_1e_1^*)$ is a self-adjoint operator. We conclude from equation (2.38) that
\begin{equation}C:=\beta_1\sqrt{-1}I \otimes e_1e_1^*+\beta_2\sqrt{-1}e_1e_1^*\otimes I+\lambda_3S(I\otimes e_1e_1^*)+\lambda_4S(e_1e_1^*\otimes I)+\lambda_5I\otimes I+\lambda_6S\end{equation} is a self-adjoint operator,
and so is the operator\begin{equation}(e_2e_2^*\otimes e_1e_1^*)C(e_2e_2^*\otimes e_1e_1^*)=\beta_1\sqrt{-1}(e_2e_2^*\otimes e_1e_1^*)+\lambda_5(e_2e_2^*\otimes e_1e_1^*).\end{equation}  Hence
\begin{equation}\beta_5=-\beta_1.\end{equation} Combining equations (2.46)-(2.47) and (2.49), we obtain that
\begin{equation}D:=\beta_2\sqrt{-1}e_1e_1^*\otimes I+\lambda_3S(I\otimes e_1e_1^*)+\lambda_4S(e_1e_1^*\otimes I)-\beta_1\sqrt{-1}I\otimes e_2e_2^*+\beta_6\sqrt{-1}S\end{equation} is a self-adjoint operator. It is easy to check that
$$(I\otimes e_2e_2^*)S(I\otimes e_2e_2^*)=S(e_2e_2^*\otimes e_2e_2^*)=e_2e_2^*\otimes e_2e_2^*.$$ Thus \begin{equation}\begin{array}{rcl}(I\otimes e_2e_2^*)D(I\otimes e_2e_2^*)&=&\beta_2\sqrt{-1}e_1e_1^*\otimes e_2e_2^*-\beta_1\sqrt{-1}I\otimes e_2e_2^*+\beta_6\sqrt{-1}e_2e_2^*\otimes e_2e_2^*\\&=&-\beta_1\sqrt{-1}e_2e_2^*\otimes e_2e_2^*+\beta_6\sqrt{-1}e_2e_2^*\otimes e_2e_2^*\end{array}\end{equation} is a self-adjoint operator, and so
\begin{equation}\beta_6=\beta_1.\end{equation}

We claim that \begin{equation} I \otimes X+X\otimes I-tr(X)I\otimes I+tr(X)S=S(X\otimes I)+S(I\otimes X)\end{equation} for all $X\in \mathcal{T(H)}.$
Indeed, since $\mathcal{H}=\bigvee\{e_1,e_2\},$ it follows that equation (2.53) holds for all $X\in \mathcal{B(H)}$  if and only if \begin{equation} I \otimes e_ie_j^*+e_ie_j^*\otimes I-S(e_ie_j^*\otimes I)-S(I\otimes e_ie_j^*)-tr(e_ie_j^*)I\otimes I+tr(e_ie_j^*)S=0\end{equation} for  $i,j=1,2.$

It is straightforward that\begin{equation}\begin{array}{rl} &I \otimes e_1e_1^*+e_1e_1^*\otimes I-S(e_1e_1^*\otimes I)-S(I\otimes e_1e_1^*)-I\otimes I+S
\\=&(I\otimes I-S)[I \otimes e_1e_1^*+e_1e_1^*\otimes I-I \otimes I]\\=&(I\otimes I-S)[e_1e_1^* \otimes e_1e_1^*-e_2e_2^*\otimes e_2e_2^*]=0.\end{array}\end{equation}
 Moreover, for all $x=t_1e_1+t_2e_2\in \mathcal{H}$ and $y=s_1e_1+s_2e_2\in \mathcal{H},$ where $t_i$ and $s_i$ $(i=1,2)$ are complex numbers, we can infer that  \begin{equation}\begin{array}{rcl}(I \otimes e_1e_2^*+e_1e_2^*\otimes I)(x\otimes y)&=&s_2x\otimes e_1+t_2e_1\otimes y\\&=&(s_2t_1+s_1t_2)e_1\otimes e_1+s_2t_2(e_1\otimes e_2+e_2\otimes e_1),\end{array}\end{equation}
and so \begin{equation}\begin{array}{rcl}&&S(I \otimes e_1e_2^*+e_1e_2^*\otimes I)(x\otimes y)\\&=&S[(s_2t_1+s_1t_2)e_1\otimes e_1+s_2t_2(e_1\otimes e_2+e_2\otimes e_1)]\\&=&(s_2t_1+s_1t_2)e_1\otimes e_1+s_2t_2(e_1\otimes e_2+e_2\otimes e_1).\end{array}\end{equation}
Then equations (2.56) and (2.57) imply $$S(I \otimes e_1e_2^*+e_1e_2^*\otimes I)=I\otimes e_1e_2^*+e_1e_2^*\otimes I,$$ which yields
$$\begin{array}{rl} I \otimes e_1e_2^*+e_1e_2^*\otimes I-S(e_1e_2^*\otimes I)-(e_1e_2^*\otimes I)S-tr(e_1e_2^*)I\otimes I+tr(e_1e_2^*)S
=0.\end{array}$$ That is to say, equation (2.54) holds for $i=1$ and $j=1,2.$
In a similar way, we can get that equation (2.54) holds for $i=2$ and $j=1,2.$ Thus (2.53) holds.

Let $\Psi: \mathcal{T}(\mathcal{H})\rightarrow\mathcal{B}(\mathcal{H}\otimes\mathcal{H})$ be defined by $$\Psi(X)=\alpha_1I \otimes X+\alpha_2X\otimes I+\alpha_3S(I\otimes X)+\alpha_4S(X\otimes I)+\alpha_5tr(X)I\otimes I+\alpha_6tr(X)S.$$ Then we get from equations (2.38), (2.46), (2.49) and (2.52)-(2.53) that \begin{equation}\begin{array}{rcl}\Phi(X)&=&\Psi(X)+\beta_3\sqrt{-1}S(I\otimes X)+\beta_4\sqrt{-1}S(X\otimes I)\\&&+\beta_1\sqrt{-1}[I \otimes X+X\otimes I-tr(X)I\otimes I+tr(X)S]\\&=&\Psi(X)+(\beta_3+\beta_1)\sqrt{-1}S(I\otimes X)+(\beta_4+\beta_1)\sqrt{-1}S(X\otimes I).\end{array}\end{equation} Put  $\mu_i=\alpha_i$ for $i=1,2,5,6,$
 $$\mu_3=\alpha_3+(\beta_3+\beta_1)\sqrt{-1} \ \ \hbox{  and }\ \ \ \mu_4=\alpha_4+(\beta_4+\beta_1)\sqrt{-1}.$$ Thus equation (2.46) implies $\mu_3=\overline{\mu_4},$ so equation (2.58) yields the desired result.

$(b)$  Let $\{e_j\}_{j=1}^n$ be
an orthonormal basis for ${\mathcal{H}}$ and $$\Phi(X)=\mu_1I \otimes X+\mu_2X\otimes I+\mu_3S(I\otimes X)+\mu_4S(X\otimes I)+\mu_5tr(X)I\otimes I+\mu_6tr(X)S,$$ where $\mu_i$ for $i=1,2,5,6$ are real numbers and complex numbers $\mu_3=\overline{\mu_4}.$

{\bf Case 1.}  Suppose that $3\leq n\leq\infty.$ If $dim\mathcal{H}=\infty$ and $P_m=\sum\limits_{i=1}^{m} e_ie_i^*,$ then $W^*$-$\lim\limits_{m\rightarrow\infty}P_m =I.$
We claim that $\Phi$ is a positive map if and only if $\Phi(e_1e_1^*)\geq0,$ which is equivalent to \begin{equation}\begin{array}{rcl}\Omega_m&:=&\mu_1P_m \otimes e_1e_1^*+\mu_2e_1e_1^*\otimes P_m+\mu_3S(P_m\otimes e_1e_1^*)\\&&+\mu_4S(e_1e_1^*\otimes P_m)+\mu_5P_m\otimes P_m+\mu_6(P_m\otimes P_m)S(P_m\otimes P_m)\geq0,\end{array}\end{equation} for $m=3,4,\cdots,n.$
Indeed, using the spectral decomposition theorem and the continuity of $\Phi$ with respect to the topologies $(\mathcal{T(H)},\|\cdot\|_1)$ and $(\mathcal{B(H\otimes H)},\|\cdot\|),$ we conclude that $\Phi$ is a positive map if and only if $\Phi(ff^*)\geq0$ for all unit vectors $f\in\mathcal{H}.$
Since $$\Phi(ff^*)=\Phi(U_fe_1e_1^*U_f^*)=(U_f\otimes U_f)\Phi(e_1e_1^*)(U_f\otimes U_f)^*,$$ where $U_f$ is a unitary operator with $U_fe_1=f,$  it follows that
$\Phi(ff^*)\geq0$ for all unit vectors $f\in\mathcal{H}$  if and only if $\Phi(e_1e_1^*)\geq0.$

Furthermore, if $dim\mathcal{H}<\infty,$ then $P_n=I,$ so $\Phi(e_1e_1^*)=\Omega_n.$
When $dim\mathcal{H}=\infty,$ $W^*$-$\lim\limits_{m\rightarrow\infty}\Omega_m =\Phi(e_1e_1^*)$ and
$$\Omega_m=(P_m\otimes P_m)\Phi(e_1e_1^*)(P_m\otimes P_m)$$ for $m=3,4,\cdots.$ Thus
$\Phi(e_1e_1^*)\geq0$ if and only if $\Omega_m\geq0$ for  $m=3,4,\cdots,n.$
It is easy to verify that \begin{equation}(P_m\otimes P_m)S(P_m\otimes P_m)=\sum\limits_{i,j=1}^{m}e_ie_j^*\otimes e_je_i^*,\end{equation}
  \begin{equation}S(P_m\otimes e_1e_1^*)=\sum\limits_{i=1}^{m}e_1 e_i^*\otimes e_i e_1^*\ \ \ \hbox{ and }\ \ S(e_1e_1^*\otimes P_m)=\sum\limits_{i=1}^{m}e_ie_1^*\otimes e_1e_i^*.\end{equation}
  Then we conclude from (2.59)-(2.61) that for  $m=3,4,\cdots,n,$
 \begin{equation}\begin{array}{rl}\Omega_m=&\mu_1\sum\limits_{i=1}^{m} e_ie_i^* \otimes e_1e_1^*+\mu_2e_1e_1^*\otimes\sum\limits_{i=1}^{m} e_ie_i^*+\mu_3\sum\limits_{i=1}^{m}e_1 e_i^*\otimes e_i e_1^*\\&+\mu_4\sum\limits_{i=1}^{m}e_i e_1^*\otimes e_1 e_i^*+\mu_5\sum\limits_{i=1}^{m} e_ie_i^*\otimes\sum\limits_{j=1}^{m} e_je_j^*+\mu_6\sum\limits_{i,j=1}^{m}e_ie_j^*\otimes e_je_i^*\\=&\sum\limits_{i=1}^{6}\mu_ie_1e_1^* \otimes e_1e_1^*+\sum\limits_{i=2}^{m} [(\mu_1+\mu_5)e_ie_i^*\otimes e_1e_1^*+(\mu_2+\mu_5)e_1e_1^*\otimes e_ie_i^*]\\+&\sum\limits_{i=2}^{m}[(\mu_3+\mu_6)e_1 e_i^*\otimes e_i e_1^*+(\mu_4+\mu_6)e_i e_1^*\otimes e_1 e_i^*]+\sum\limits_{i,j=2}^{m}(\mu_5 e_ie_i^*\otimes e_je_j^*+\mu_6e_ie_j^*\otimes e_je_i^*).\end{array}\end{equation}
 For simplicity, we denote by $\Pi_1:=\sum\limits_{i=1}^{6}\mu_ie_1e_1^* \otimes e_1e_1^*,$ $$\Pi_{m+1}:=\sum\limits_{i,j=2}^{m}(\mu_5 e_ie_i^*\otimes e_je_j^*+\mu_6e_ie_j^*\otimes e_je_i^*)$$ and for $i=2,3,\cdots,m,$ $$\begin{array}{rcl}\Pi_i&:=&(\mu_1+\mu_5)e_ie_i^*\otimes e_1e_1^*\\&&+(\mu_2+\mu_5)e_1e_1^*\otimes e_ie_i^*+(\mu_3+\mu_6)e_1 e_i^*\otimes e_i e_1^*+(\mu_4+\mu_6)e_i e_1^*\otimes e_1 e_i^*.\end{array}$$ It is easy to check that
  $$ \Omega_m=\sum\limits_{i=1}^{m+1} \Pi_i\ \ \ \hbox{ and } \ \ \ \Pi_k\Pi_l=0=\Pi_l\Pi_k \ \ \hbox{ for }\ k\neq l=1,2,\cdots m+1.$$ Thus $\Omega_m\geq0$ if and only if $\Pi_i\geq0$ for $i=1,2,\cdots m+1.$
  Since $m\geq 3$ and the eigenvalues of $\sum\limits_{i,j=2}^{m}e_ie_j^*\otimes e_je_i^*$ acting on the subspace $ran(\sum\limits_{i,j=2}^{m}e_ie_i^*\otimes e_je_j^* )$ are $\{1,-1\},$
  it follows that the eigenvalues of $\Pi_{m+1}$ acting on the subspace $ran(\sum\limits_{i,j=2}^{m}e_ie_i^*\otimes e_je_j^* )$ are $\{\mu_5+\mu_6,\mu_5-\mu_6\}.$
  So $\Pi_{m+1}\geq0$ if and only if $\mu_5\geq |\mu_6|.$
   Clearly, $\Pi_1\geq0$ if and only if $\sum_{i=1}^{6}\mu_i\geq0.$ Moreover, it is easy to verify that $\Pi_i\geq0$ for $i=2,3,\cdots m$ if and only if
 $$\left(\begin{array}{cc}
\mu_1+\mu_5&\overline{\mu_3}+\mu_6\\ \mu_3+\mu_6&\mu_2+\mu_5\end{array}\right)\geq0.$$

{\bf Case 2.} Let $dim \mathcal{H}=2$ and $\{e_j\}_{j=1}^2$ is
an orthonormal basis for ${\mathcal{H}}.$ As shown above, $\Phi$ is a positive map if and only if $\Phi(e_1e_1^*)\geq0.$ Then a direct calculation implies
$$\begin{array}{rcl}\Phi(e_1e_1^*)&=&\sum\limits_{i=1}^{6}\mu_ie_1e_1^* \otimes e_1e_1^*+(\mu_1+\mu_5)e_2e_2^*\otimes e_1e_1^*+(\mu_2+\mu_5)e_1e_1^*\otimes e_2e_2^*\\&+&(\mu_6+\mu_5)e_2e_2^*\otimes e_2e_2^*+(\mu_3+\mu_6)e_1 e_2^*\otimes e_2 e_1^*+(\overline{\mu_3}+\mu_6)e_2 e_1^*\otimes e_1 e_2^*\\&\simeq&\left(\begin{array}{cccc}
\sum\limits_{i=1}^{6}\mu_i&0&0&0\\ 0&\mu_1+\mu_5&\overline{\mu_3}+\mu_6&0
\\0&\mu_3+\mu_6&\mu_2+\mu_5&0\\ 0&0&0&\mu_6+\mu_5 \end{array}\right).\end{array}$$ Thus $\Phi(e_1e_1^*)\geq0 $ if and only if $\mu_5+\mu_6\geq0,$ $\sum_{i=1}^{6}\mu_i\geq0$  and
 $$\left(\begin{array}{cc}
\mu_1+\mu_5&\overline{\mu_3}+\mu_6 \\ \mu_3+\mu_6 &\mu_2+\mu_5\end{array}\right)\geq0.$$ $\Box$

{\bf Corollary 3.} Let $dim\mathcal{H}=n<\infty$ and  $\Phi:\mathcal{B(H)}\rightarrow \mathcal{B(H\otimes H)}$ be a linear map.

$(a)$ $\Phi$ is covariance under unitary evolution and satisfies the broadcasting
condition (1.1)  if and only if  there exist complex numbers $\lambda_i$ for $i=1,2,\cdots,6$ with $\lambda_1=\lambda_2,$ $n\lambda_1+\lambda_3+\lambda_4=1$ and $n\lambda_5+\lambda_2+\lambda_6=0$ such that $$\Phi(X)=\lambda_1I \otimes X+\lambda_2X\otimes I+\lambda_3S(I\otimes X)+\lambda_4S(X\otimes I)+\lambda_5tr(X)I\otimes I+\lambda_6tr(X)S,$$ where $S\in\mathcal{B(H\otimes H)}$ is the swap operator.

$(b)$ If $\Phi$ is covariance under unitary evolution and satisfies the broadcasting
condition (1.1), then $\Phi$ is not a positive map.

{\bf Proof.} $(a)$ It is easy to verify that $tr_1(S)=tr_2(S)=I$ and for all $X\in\mathcal{B(H)},$   $$tr_1[S(I\otimes X)]=tr_2[S(I\otimes X)]=tr_1[S(X\otimes I)]=tr_2[S(X\otimes I)]=X.$$
Then $(a)$ follows from Theorem 1 and a direct calculation.

$(b)$ Suppose that $\Phi$ is a positive map.  If $n\geq 3,$ then we conclude from above $(a)$ and Theorem 2 that $\lambda_5+\lambda_2\geq0$ and $\lambda_5+\lambda_6\geq0.$ So $n\lambda_5+\lambda_2+\lambda_6=0$ implies that $$\lambda_5=\lambda_2=\lambda_6=0.$$ Moreover, $\lambda_4=\overline{\lambda_3}=0$ follows from
the inequality  $$\left(\begin{array}{cc}
\lambda_1+\lambda_5&\overline{\lambda_3}+\lambda_6\\ \lambda_3+\lambda_6&\lambda_2+\lambda_5\end{array}\right)\geq0.$$
Thus $(a)$ yields that $$0=n\lambda_1+\lambda_3+\lambda_4=1,$$ which is a contradiction.

If $n=2,$ then $2\lambda_5+\lambda_2+\lambda_6=0.$ It entails that  $$\lambda_5+\lambda_2=\lambda_5+\lambda_6=0,$$ and so $$\lambda_1=\lambda_2=\lambda_6=-\lambda_5 \ \hbox{ and }\ \
   \lambda_3+\lambda_6=0.$$ Thus $\lambda_4=\lambda_3=\lambda_5,$ and hence $$0=2\lambda_1+\lambda_3+\lambda_4=1.$$ It is a contradiction. $\Box$

{\bf Remark 4.} Let $\Phi:\mathcal{T(H)}\rightarrow \mathcal{B(H\otimes H)}$ be a linear map and continuous with respect to the topologies $(\mathcal{T(H)},\|\cdot\|_1)$ and $(\mathcal{B(H\otimes H)},\|\cdot\|).$ Suppose that $\Phi$ is covariance under unitary evolution. Then Theorem 1 implies that there exist complex numbers $\lambda_i$ for $i=1,2,\cdots,6$ such that $$\Phi(X)=\lambda_1I \otimes X+\lambda_2X\otimes I+\lambda_3S(I\otimes X)+\lambda_4S(X\otimes I)+\lambda_5tr(X)I\otimes I+\lambda_6tr(X)S,$$ where $S\in\mathcal{B(H\otimes H)}$ is the swap operator. We claim that if $dim \mathcal{H}\geq 3,$ then the complex numbers $\lambda_i$ are unique.
Indeed, if $e_1,$ $e_2$ and $e_3$ are orthogonal unit vectors of $\mathcal{H}$ and if there exist other complex numbers $\mu_i$ for $i=1,2,\cdots,6$ such that $$\Phi(X)=\mu_1I \otimes X+\mu_2X\otimes I+\mu_3S(I\otimes X)+\mu_4S(X\otimes I)+\mu_5tr(X)I\otimes I+\mu_6tr(X)S,$$
then $$(\mu_1-\lambda_1)I \otimes e_1e_2^*+(\mu_2-\lambda_2)e_1e_2^*\otimes I+(\mu_3-\lambda_3)S(I\otimes e_1e_2^*)+(\mu_4-\lambda_4)S(e_1e_2^*\otimes I)=0,$$ so $$(\mu_1-\lambda_1)e_3e_3^* \otimes e_1e_1^*+ (\mu_3-\lambda_3)S(e_3e_3^*\otimes e_1e_1^*)=0,$$ which yields $\mu_1=\lambda_1$ and $\mu_3=\lambda_3.$
In a similar way, we get that $\mu_2=\lambda_2$ and $\mu_4=\lambda_4.$ Thus $$(\mu_5-\lambda_5) I\otimes I+(\mu_6-\lambda_6)S=0,$$ which implies that $\mu_5=\lambda_5$ and $\mu_6=\lambda_6.$ However, if $dim \mathcal{H}=2,$ then we conclude from equation (2.53) that $$\begin{array}{rc}&\Phi(X)=(\lambda_1+\mu)I \otimes X+(\lambda_2+\mu)X\otimes I+(\lambda_3-\mu)S(I\otimes X)+(\lambda_4-\mu)S(X\otimes I)\\&+(\lambda_5-\mu)tr(X)I\otimes I+(\lambda_6+\mu)tr(X)S\end{array}$$ for any complex number $\mu.$  $\Box$

In the following, we get the specific structure of the linear map $\Phi$ that is covariance under unitary evolution and invariance under permutation.

{\bf Corollary 5.} Let $\Phi:\mathcal{T(H)}\rightarrow \mathcal{B(H\otimes H)}$ be a linear map and continuous with respect to the topologies $(\mathcal{T(H)},\|\cdot\|_1)$ and $(\mathcal{B(H\otimes H)},\|\cdot\|).$ Then $\Phi$ is covariance under unitary evolution and invariance under permutation if and only if there exist complex numbers $\mu_i$ for $i=1,2,3,4$ such that $$\Phi(X)=\mu_1(I \otimes X+X\otimes I)+\mu_2[S(I\otimes X)+S(X\otimes I)]+\mu_3tr(X)I\otimes I+\mu_4tr(X)S,$$ where $S\in\mathcal{B(H\otimes H)}$ is the swap operator.

{\bf Proof.}  Sufficiency is obvious.

Necessity. Using Theorem 1, we conclude that there exist complex numbers $\lambda_i$ for $i=1,2,\cdots,6$ such that $$\Phi(X)=\lambda_1I \otimes X+\lambda_2X\otimes I+\lambda_3S(I\otimes X)+\lambda_4S(X\otimes I)+\lambda_5tr(X)I\otimes I+\lambda_6tr(X)S.$$ Let $e_1$ and $e_2$ be unit vectors with $e_1\bot e_2.$ Since $\Phi$ is invariance under permutation, it follows that $\Phi(e_1e_1^*)=S\Phi(e_1e_1^*)S,$ so $$\begin{array}{rl}&\lambda_1I \otimes e_1e_1^*+\lambda_2e_1e_1^*\otimes I+\lambda_3S(I\otimes e_1e_1^*)+\lambda_4S(e_1e_1^*\otimes I)+\lambda_5I\otimes I+\lambda_6S\\=&\lambda_1 e_1e_1^* \otimes I+\lambda_2 I\otimes e_1e_1^*+\lambda_3S(e_1e_1^*\otimes I)+\lambda_4S(I\otimes e_1e_1^*)+\lambda_5I\otimes I+\lambda_6S,\end{array}$$ which implies $$\begin{array}{rl}(\lambda_1-\lambda_2)(e_1e_1^*\otimes I -I \otimes e_1e_1^*)+(\lambda_3-\lambda_4)S(e_1e_1^*\otimes I -I\otimes e_1e_1^*)=0.\end{array}$$
Thus $$\begin{array}{rl}&(\lambda_1-\lambda_2)e_1\otimes e_2+(\lambda_3-\lambda_4)e_2\otimes e_1\\=&[(\lambda_1-\lambda_2)(e_1e_1^*\otimes I -I \otimes e_1e_1^*)+(\lambda_3-\lambda_4)S(e_1e_1^*\otimes I -I\otimes e_1e_1^*)]e_1\otimes e_2=0,\end{array}$$ so $\lambda_1=\lambda_2$ and $\lambda_3=\lambda_4.$ Then
we denote by $$\mu_1=\lambda_1, \  \ \ \mu_2=\lambda_3, \ \ \ \ \mu_3=\lambda_5\ \ \hbox{ and }\ \mu_4=\lambda_6$$ and get the desired result. $\Box$

Analogously, we give a concrete description of the linear maps $\Phi,$ which are covariance under unitary evolution and consistency with classical broadcasting.

{\bf Corollary 6.} Let $\Phi:(\mathcal{T}(\mathcal{H}), \tau_{1})\rightarrow (\mathcal{B}(\mathcal{H}\otimes\mathcal{H}), \tau)$ be a continuous linear  map. Then $\Phi$ is covariance under unitary evolution and consistency with classical broadcasting if and only if there exists complex numbers $\mu$ such that $$\Phi(X)=\mu S(I\otimes X)+(1-\mu)S(X\otimes I),$$ where $S\in\mathcal{B(H\otimes H)}$ is the swap operator.

{\bf Proof.}  Sufficiency is easy to verify.

Necessity.  Using Theorem 1, we get that there exist complex numbers $\lambda_i$ for $i=1,2,\cdots,6$ such that $$\Phi(X)=\lambda_1I \otimes X+\lambda_2X\otimes I+\lambda_3S(I\otimes X)+\lambda_4S(X\otimes I)+\lambda_5tr(X)I\otimes I+\lambda_6tr(X)S.$$ Let $\{e_j\}_{j=1}^n$ be
a fixed orthonormal basis for ${\mathcal{H}},$ where $2\leq n\leq\infty.$ Since $\Phi$ is consistency with classical broadcasting, it follows that \begin{equation}\begin{array}{rl}&\lambda_1I \otimes e_1e_1^*+\lambda_2e_1e_1^*\otimes I+\lambda_3e_1e_1^*\otimes e_1e_1^*+\lambda_4e_1e_1^*\otimes e_1e_1^*+\lambda_5I\otimes I+\lambda_6\sum\limits_{i=1}^{s}e_ie_i^*\otimes e_ie_i^*\\=&(\mathcal{D}\otimes\mathcal{D})\circ\Phi\circ\mathcal{D}(e_1e_1^*)
=\mathcal{B}_{cl}(e_1e_1^*)=e_1e_1^*\otimes e_1e_1^*.\end{array}\end{equation}

{\bf Case 1.} If $n\geq 3,$ then equation (2.63) implies $$\begin{array}{rl}&\lambda_5e_3e_3^*\otimes e_2e_2^*\\=&[\lambda_1I \otimes e_1e_1^*+\lambda_2e_1e_1^*\otimes I+\lambda_3e_1e_1^*\otimes e_1e_1^*+\lambda_4e_1e_1^*\otimes e_1e_1^* +\lambda_5I\otimes I\\&+\lambda_6\sum\limits_{i=1}^{s}e_ie_i^*\otimes e_ie_i^*]e_3e_3^*\otimes e_2e_2^*=0.\end{array}$$ So $\lambda_5=0.$ In a similar way, we also get that
$$\lambda_1=\lambda_2=\lambda_6=0.$$ Thus $\lambda_3+\lambda_4=1$ follows from equation (2.63).
Denoting by $\mu=\lambda_3,$ we get the desired equation.

{\bf Case 2.}  If $n=2,$ then equation (2.63) yields that $$\begin{array}{rl}&(\sum\limits_{i=1}^{6}\lambda_i-1)e_1e_1^*\otimes e_1e_1^*+(\lambda_1+\lambda_5)e_2e_2^*\otimes e_1e_1^*+(\lambda_2+\lambda_5)e_1e_1^*\otimes e_2e_2^*\\ +&(\lambda_5+\lambda_6)e_2e_2^*\otimes e_2e_2^*=0.\end{array}$$ Thus
$$\lambda_1=\lambda_2=\lambda_6=-\lambda_5\ \ \ \hbox{and} \ \ 2\lambda_1+\lambda_3+\lambda_4=1,$$ which implies that
  \begin{equation}\Phi(X)=\lambda_1[I \otimes X+X\otimes I-tr(X)I\otimes I +tr(X)S]+\lambda_3 S(I\otimes X)+\lambda_4S(X\otimes I).\end{equation} Since $n=2,$ it follows from equation (2.53) that $$I \otimes X+X\otimes I-tr(X)I\otimes I+tr(X)S=S(X\otimes I)+S(I\otimes X)$$ for all $X\in \mathcal{T(H)}.$
  Then equation (2.64) becomes  $$\Phi(X)=(\lambda_3+\lambda_1) S(I\otimes X)+(1-\lambda_3-\lambda_1) S(X\otimes I).$$
   Setting $\mu=\lambda_1+\lambda_3,$ we conclude that the desired equation holds.   $\Box$

Combining Corollaries 5 and 6, we can obtain the following result.

{\bf Corollary 7.} Let $\Phi:(\mathcal{T}(\mathcal{H}), \tau_{1})\rightarrow (\mathcal{B}(\mathcal{H}\otimes\mathcal{H}), \tau)$ be a continuous linear map.
Then $\Phi$ is covariance under unitary evolution, invariance under permutation and consistency with classical broadcasting if and only if $\Phi(X)=\frac{1}{2}[S(I\otimes X)+S(X\otimes I)],$ where $S\in\mathcal{B(H\otimes H)}$ is the swap operator.

{\bf Remark 8.} As is well known if $\mathcal{H}$ is finite dimensional, then $\mathcal{T}(\mathcal{H})=\mathcal{B}(\mathcal{H})$ and the linear map $\Phi:\mathcal{T(H)}\rightarrow \mathcal{B(H\otimes H)}$ is (automatic) continuous. 
 Thus Corollary 7 is an extension of [14, Theorem 1] to a separable Hilbert space $\mathcal{H}.$ Particularly,
 Corollary 6 shows that for a finite dimensional Hilbert space $\mathcal{H},$ if $\Phi$ is covariance under unitary evolution and consistency with classical broadcasting, then
  $$ tr_1[\Phi(X)]=tr_2[\Phi(X)]=X  \ \ \ \hbox{ for all }\ \ \ X\in \mathcal{T(H)}.$$ That is, the broadcasting condition (1.1) holds for $\Phi.$ Note that the broadcasting condition (1.1) is an explicit assumption in [14, Theorem 1].  $\Box$

\section{ Some extensions of Theorem 1 }

Let $\mathcal{K(H)}$ be the $C^*$-algebra of all compact operators on $\mathcal{H}.$  Based on Theorem 1, we can also characterize the structure of the norm continuous linear maps $\Psi:\mathcal{K(H)}\rightarrow \mathcal{B(H\otimes H)},$ which are covariance under unitary evolution.

{\bf Proposition 9.} Let ${\mathcal{H}}$ be an infinite-dimensional separable Hilbert space and $\Psi:\mathcal{K(H)}\rightarrow \mathcal{B(H\otimes H)}$ be a norm continuous linear map. Then the following statements are equivalent

$(a)$ $\Psi$ is covariance under unitary evolution: $\Psi(UXU^*)=(U\otimes U)\Psi(X)(U^*\otimes U^*)$ for all $X\in\mathcal{K(H)}$ and $U\in\mathcal{U(H)};$

$(b)$ $\Psi(Uee^*U^*)=(U\otimes U)\Psi(ee^*)(U^*\otimes U^*)$ for a unit vector $e\in\mathcal{H}$ and all $U\in\mathcal{U(H)};$

$(c)$ There exist complex numbers $\lambda_i$ for $i=1,2,3,4$ such that $$\Psi(X)=\lambda_1I \otimes X+\lambda_2X\otimes I+\lambda_3S(I\otimes X)+\lambda_4S(X\otimes I),$$ where $S\in\mathcal{B(H\otimes H)}$ is the swap operator.

{\bf Proof.} $(a)\Longrightarrow(b)$ and $(c)\Longrightarrow(a)$ are clear.

$(b)\Longrightarrow(c).$ Let $\{e_j\}_{j=1}^\infty$ be
an orthonormal basis for ${\mathcal{H}}.$ We conclude from the proof of Theorem 1 that there exist complex numbers $\lambda_1,\lambda_2,\cdots,\lambda_6$ such that \begin{equation}\Psi(e_je_j^*)=\lambda_1I \otimes e_je_j^*+\lambda_2e_je_j^*\otimes I+\lambda_3S(I\otimes e_je_j^*)+\lambda_4S(e_je_j^*\otimes I)+\lambda_5I\otimes I+\lambda_6S,\end{equation}
for $j=1,2,3,\cdots.$ Setting $A_n=\sum_{i=1}^{n} \frac{1}{i}e_ie_i^*,$ we get from equation (3.1) that
 \begin{equation}\Psi(A_n)=\lambda_1I \otimes A_n+\lambda_2A_n\otimes I+\lambda_3S(I\otimes A_n)+\lambda_4S(A_n\otimes I)+\lambda_5\sum_{i=1}^{n} \frac{1}{i}I\otimes I+\lambda_6\sum_{i=1}^{n} \frac{1}{i}S.\end{equation}  Clearly, $\lim\limits_{n\rightarrow\infty}\|A_n-A\|=0,$ where $$A=\sum_{i=1}^{\infty} \frac{1}{i}e_ie_i^*\in\mathcal{K(H)}.$$ Since $\Psi$ is norm continuous, we have $\lim\limits_{n\rightarrow\infty}\|\Psi(A_n)-\Psi(A)\|=0,$  so equation (3.2) yields \begin{equation}\lim\limits_{n\rightarrow\infty}\| \lambda_5\sum_{i=1}^{n} \frac{1}{i}I\otimes I+\lambda_6\sum_{i=1}^{n} \frac{1}{i}S-B\|=0,\end{equation} where  $$B:=\Psi(A)-[\lambda_1I \otimes A+\lambda_2A\otimes I+\lambda_3S(I\otimes A)+\lambda_4S(A\otimes I)].$$
 Then equation (3.3) implies  \begin{equation}\lim\limits_{n\rightarrow\infty}\|\sum_{i=1}^{n} \frac{1}{i} (\lambda_5+\lambda_6)(I\otimes I+S)-(I\otimes I+S)B\|=0 \end{equation}
and \begin{equation}\lim\limits_{n\rightarrow\infty}\|\sum_{i=1}^{n} \frac{1}{i} (\lambda_5-\lambda_6)(I\otimes I-S)-(I\otimes I-S)B\|=0.\end{equation} Thus \begin{equation}2|\lambda_5+\lambda_6|\lim\limits_{n\rightarrow\infty}\sum_{i=1}^{n} \frac{1}{i}=\lim\limits_{n\rightarrow\infty}\|\sum_{i=1}^{n} \frac{1}{i} (\lambda_5+\lambda_6)(I\otimes I+S)\|=\|(I\otimes I+S)B\| \end{equation}
and \begin{equation}2|\lambda_5-\lambda_6|\lim\limits_{n\rightarrow\infty}\sum_{i=1}^{n} \frac{1}{i}=\lim\limits_{n\rightarrow\infty}\sum_{i=1}^{n} \frac{1}{i}\| (\lambda_5-\lambda_6)(I\otimes I-S)\|= \|(I\otimes I-S)B\|.\end{equation}
Then $\lambda_5=\lambda_6=0$ follows from equations (3.6) and (3.7),
 and by  equation (3.1), we get that
 $$\Psi(ee^*)=\lambda_1I \otimes ee^*+\lambda_2ee^*\otimes I+\lambda_3S(I\otimes ee^*)+\lambda_4S(ee^*\otimes I)$$ for all unit vectors $e\in\mathcal{H}.$ The remaining proof is similar to that of Theorem 1. $\Box$

{\bf Proof of Theorem 3.} $(a)\Longrightarrow(b)$ and $(c)\Longrightarrow(a)$ are obvious.

$(b)\Longrightarrow(c).$ Let $A_n,$ $A$ and $B$ be as in the proof of Proposition 9. Since $\Psi$ is a $W^*$ continuous linear map, it follows in a similar way to equation (3.3) that \begin{equation}W^*\hbox{-}\lim\limits_{n\rightarrow\infty} (\lambda_5\sum_{i=1}^{n} \frac{1}{i}I\otimes I+\lambda_6\sum_{i=1}^{n} \frac{1}{i}S)=B.\end{equation}
So  \begin{equation}W^*\hbox{-}\lim\limits_{n\rightarrow\infty}\sum_{i=1}^{n} \frac{1}{i} (\lambda_5+\lambda_6)(I\otimes I+S)=(I\otimes I+S)B \end{equation}
and \begin{equation}W^*\hbox{-}\lim\limits_{n\rightarrow\infty}\sum_{i=1}^{n} \frac{1}{i} (\lambda_5-\lambda_6)(I\otimes I-S)=(I\otimes I-S)B.\end{equation}
It is clear that $\frac{I\otimes I+S}{2}$ and $\frac{I\otimes I-S}{2}$ are non-zero orthogonal projections.
Suppose that $f,g\in\mathcal{H}$ are unit vectors such that $\frac{I\otimes I+S}{2}f=f$ and $\frac{I\otimes I-S}{2}g=g.$ Then equations (3.9) and (3.10) imply that \begin{equation}2(\lambda_5+\lambda_6)\lim\limits_{n\rightarrow\infty}\sum_{i=1}^{n} \frac{1}{i}=\lim\limits_{n\rightarrow\infty}\sum_{i=1}^{n} \frac{1}{i} (\lambda_5+\lambda_6)tr[(I\otimes I+S)ff^*]=tr[(I\otimes I+S)Bff^*] \end{equation}
and \begin{equation}2(\lambda_5-\lambda_6)\lim\limits_{n\rightarrow\infty}\sum_{i=1}^{n} \frac{1}{i}=\lim\limits_{n\rightarrow\infty}\sum_{i=1}^{n} \frac{1}{i}(\lambda_5-\lambda_6)tr[(I\otimes I-S)gg^*]= tr[(I\otimes I-S)Bgg^*].\end{equation}
Thus  $\lambda_5=\lambda_6=0,$ and hence  for all unit vectors $e\in\mathcal{H},$   $$\Psi(ee^*)=\lambda_1I \otimes ee^*+\lambda_2ee^*\otimes I+\lambda_3S(I\otimes ee^*)+\lambda_4S(ee^*\otimes I).$$ Therefore,
\begin{equation}\Psi(X)=\lambda_1I \otimes X+\lambda_2X\otimes I+\lambda_3S(I\otimes X)+\lambda_4S(X\otimes I),\end{equation} for all finite rank operators $X.$ It is known that the set of all finite rank operators is dense in $\mathcal{B(H)}$ with respect to $W^*$-topology, so equation (3.13) implies that $(c)$ holds as desired. $\Box$

 Combining the proof of Corollary 5, Corollary 6 and Theorem 3, we get another extension for the virtual broadcasting map.

{\bf Corollary 10.} Let ${\mathcal{H}}$ be an infinite-dimensional separable Hilbert space and $\Psi:\mathcal{B(H)}\rightarrow \mathcal{B(H\otimes H)}$ be a continuous linear map in the $W^*$-topology. Suppose that $S\in\mathcal{B(H\otimes H)}$ is the swap operator.

$(a)$ $\Psi$ is covariance under unitary evolution and consistency with classical broadcasting if and only if   there exists complex numbers $\mu$ such that $\Psi(X)=\mu S(I\otimes X)+(1-\mu)S(X\otimes I).$

$(b)$  $\Psi$ is covariance under unitary evolution, invariance under permutation and consistency with classical broadcasting if and only if
$\Psi(X)=\frac{ 1}{2}[S(I\otimes X)+S(X\otimes I)].$

\section{Completely positivity and completely bounded norm of $\Psi$ }

Let $\Phi:$ ${\mathcal{B(H)}}\longrightarrow
  {\mathcal{B(K)}}$ be a linear map. For an integer $n\geq1,$ suppose that $\mathcal{M}_n({\mathcal{B(H)}})$ is the von Neumann algebra of all
$n\times n$ matrices
  whose entries are operators in ${\mathcal{B(H)}}$ and $\Phi_n:$ $\mathcal{M}_n({\mathcal{B(H)}})\longrightarrow
\mathcal{M}_n({\mathcal{B(K)}})$ is the linear map defined by 
 $$\Phi_n([A_{i,j}]_{n, n})=[\Phi(A_{i,j})]_{n, n}, \ \ \ \hbox {   for  }
 [A_{i,j}]_{n, n}\in\mathcal{M}_n({\mathcal{B(H)}}).$$ $\Phi$
is said to be $n$ positive if $\Phi_n([A_{i,j}]_{n, n})\geq0,$ when $[A_{i,j}]_{n, n}\geq0.$
 If $\Phi_n$ is positive for all $n\geq 1,$ then $\Phi$ is called completely positive. Moreover,
 we call $\Phi$ completely bounded if $\sup_{n}\|\Phi_n\|<\infty,$ where for $n=1,2,\cdots,$  $$\|\Phi_n\|=sup\{\|\Phi_n([A_{i,j}]_{n, n})\|: [A_{i,j}]_{n, n}\in\mathcal{M}_n({\mathcal{B(H)}}) \hbox{ with } \|[A_{i,j}]_{n, n}\|\leq1\}.$$ In this case, denote by $\|\Phi\|_{cb}=\sup_{n}\|\Phi_n\|$ the completely bounded norm of $\Psi$ ([16]).
Suppose that $\Psi:\mathcal{B(H)}\rightarrow \mathcal{B(H\otimes H)}$ is a continuous linear map in the $W^*$-topology and covariance under unitary evolution. In the following, we show that $\Psi$ is a positive map if and only if $\Psi$ is a completely positive map.

{\bf Corollary 11.} Let $\Psi:\mathcal{B(H)}\rightarrow \mathcal{B(H\otimes H)}$ be a linear map with the form  $$\Psi(X)=\lambda_1I \otimes X+\lambda_2X\otimes I+\lambda_3S(I\otimes X)+\lambda_4S(X\otimes I),$$ where $S\in\mathcal{B(H\otimes H)}$ is the swap operator and $\lambda_i$ are complex numbers for $i=1,2,3,4.$  Then the following statements are equivalent:

$(a)$  $\Psi$ is a positive map;

$(b)$ $\lambda_1\geq0,$ $\lambda_2\geq0,$ $\lambda_4=\overline{\lambda_3}$ and $\lambda_1\lambda_2\geq |\lambda_3|^2;$

$(c)$  $\Psi$ is a completely positive map.

{\bf Proof.} $(c)\Longrightarrow(a)$ is obvious.  $(a)\Longrightarrow(b).$  Using Corollary 3, we get $\lambda_4=\overline{\lambda_3}.$ Then Theorem 2 implies  that  $\left(\begin{array}{cc}
\lambda_1 &\overline{\lambda_3 }\\  \lambda_3 &\lambda_2\end{array}\right)\geq0,$ so $\lambda_1\geq0,$ $\lambda_2\geq0$  and $\lambda_1\lambda_2\geq |\lambda_3|^2.$

$(b)\Longrightarrow(c).$ If $\lambda_1=0,$ then inequality $\lambda_1\lambda_2\geq |\lambda_3|^2$ yields $\lambda_4=\overline{\lambda_3}=0.$
Since $\lambda_2\geq0,$ it follows that $\Psi(X)=\lambda_2X\otimes I$  is a completely positive map.
Suppose that $\lambda_1>0$  and $\mu=\frac{|\lambda_3|^2}{\lambda_1}.$ Then $\lambda_2\geq\mu$
and $$\begin{array}{rcl}\Psi(X)&=&(\lambda_2-\mu)X\otimes I+\lambda_1I \otimes X+\mu X\otimes I+\lambda_3S(I\otimes X)+\lambda_4S(X\otimes I)\\&=&(\lambda_2-\mu)X\otimes I+(\sqrt{\lambda_1}\overline{\omega}I \otimes I,\sqrt{\mu}S)\left(\begin{array}{cc}
  I\otimes X&I\otimes X \\  I\otimes X & I\otimes X\end{array}\right)\left(\begin{array}{c}
\sqrt{\lambda_1}\omega I \otimes I \\  \sqrt{\mu}  S\end{array}\right),\end{array}$$ where $\lambda_3=|\lambda_3|\omega.$ Clearly, the linear map $X\longrightarrow\left(\begin{array}{cc}
    X&  X \\   X &   X\end{array}\right)$ is completely positive from $\mathcal{B(H)}$ into $\mathcal{B( H\oplus H)},$ which implies that $$\Lambda(X)=(\sqrt{\lambda_1}\overline{\omega}I \otimes I,\sqrt{\mu}S)\left(\begin{array}{cc}
  I\otimes X&I\otimes X \\  I\otimes X & I\otimes X\end{array}\right)\left(\begin{array}{c}
\sqrt{\lambda_1}\omega I \otimes I \\  \sqrt{\mu}  S\end{array}\right)$$ is completely positive, so
 $\Psi$ is a completely positive map.
$\Box$

Combining Theorem 3 and proof of Corollary 5, we get the following result.

{\bf Corollary 12.} Let ${\mathcal{H}}$ be an infinite-dimensional Hilbert space and $\Psi:\mathcal{B(H)}\rightarrow \mathcal{B(H\otimes H)}$ be a continuous linear map in the $W^*$-topology.
Suppose $\Psi$ is covariance under unitary evolution and invariance under permutation. Then

$(a)$ $\Psi$ is unital   $(\Psi(I)=I\otimes I)$
if and only if $\Psi(X)=\frac{ 1}{2}(I\otimes X+X\otimes I).$

$(b)$ $\Psi(I)=S$
if and only if $\Psi(X)=\frac{ 1}{2}[S(I\otimes X)+S(X\otimes I)],$ where $S\in\mathcal{B(H\otimes H)}$ is the swap operator.

$(c)$ $\|\Psi\|_{cb}=\|\Psi\|=\|\Psi(I)\|.$

{\bf Proof.}  Since $\Psi$ is covariance under unitary evolution and invariance under permutation, it follows from Theorem 3 and proof of Corollary 5 that there exist complex numbers $\mu_1$ and $\mu_2$ such that \begin{equation}\Psi(X)=\mu_1(I \otimes X+X\otimes I)+\mu_2[S(I\otimes X)+S(X\otimes I)].\end{equation} Then it can be readily verified that $(a)$ and $(b)$ are in force.

$(c).$ We denote by $Q:=\frac{I\otimes I+S}{2}.$ Then $Q$ is non-zero orthogonal projections and $Q^\perp=I\otimes I-Q=\frac{I\otimes I-S}{2}.$  It is easy to see that \begin{equation} Q(I \otimes X)Q=\frac{1}{4}[I \otimes X+X\otimes I+S(I\otimes X)+S(X\otimes I)] \end{equation} and \begin{equation}   Q^\perp(I\otimes X)Q^\perp=\frac{1}{4}[I \otimes X+X\otimes I-S(I\otimes X)-S(X\otimes I)].\end{equation}
Combining equations (4.1)-(4.3), we get that \begin{equation}\Psi(X)=2(\mu_1+\mu_2)Q(I\otimes X)Q+2(\mu_1-\mu_2)Q^\perp(I\otimes X)Q^\perp.\end{equation}
For any $[A_{i,j}]_{n, n}\in\mathcal{M}_n({\mathcal{B(H)}}),$  we conclude from equation (4.4) that  $$\begin{array}{rcl}\|\Psi_n([A_{i,j}]_{n, n})\|&=&\|2(\mu_1+\mu_2)[Q(I\otimes A_{i,j})Q+2(\mu_1-\mu_2)Q^\perp(I\otimes A_{i,j})Q^\perp]_{n, n}\|\\&=&\|[2(\mu_1+\mu_2)Q(I\otimes A_{i,j})Q]_{n, n}+2(\mu_1-\mu_2)[Q^\perp(I\otimes A_{i,j})Q^\perp]_{n, n}\|\\&=&max\{2|\mu_1+\mu_2|\|[Q(I\otimes A_{i,j})Q]_{n, n}\|,\ 2|\mu_1-\mu_2|\|[Q^\perp(I\otimes A_{i,j})Q^\perp]_{n, n}\|\}\\&\leq&max\{2|\mu_1+\mu_2|\|[I\otimes A_{i,j}]_{n, n}\|,\ 2|\mu_1-\mu_2|\|[I\otimes A_{i,j}]_{n, n}\|\}\\&=&max\{2|\mu_1+\mu_2|,\ 2|\mu_1-\mu_2|\}\|[A_{i,j}]_{n, n}\|.\end{array} $$
Thus for $n=1,2,\cdots,$ $$\|\Psi_n \|\leq max\{2|\mu_1+\mu_2|,\ 2|\mu_1-\mu_2|\},$$
and so  \begin{equation}\|\Psi\|_{cb}\leq max\{2|\mu_1+\mu_2|,\ 2|\mu_1-\mu_2|\}.\end{equation}
Moreover, equation (4.1) implies \begin{equation}\|\Psi\|_{cb}\geq\|\Psi\|\geq\|\Psi(I)\|= \|2\mu_1(I \otimes I)+2\mu_2S\|=max\{2|\mu_1+\mu_2|,\ 2|\mu_1-\mu_2|\}.\end{equation}
Then $\|\Psi\|_{cb}=\|\Psi\|=\|\Psi(I)\|$ follows from equations (4.5) and (4.6).   $\Box$

More generally, we consider how to characterize the upper bounded  of $\|\Psi\|_{cb},$ when $\Psi:\mathcal{B(H)}\rightarrow \mathcal{B(H\otimes H)}$ is a continuous linear map in the $W^*$-topology and covariance under unitary evolution.  A technical lemma is needed.

{\bf Lemma 13.} Let $P\in\mathcal{B(H)}$ be an orthogonal projection and $\mu_i$ be complex numbers for $i=1,2,3,4.$ If $A\in\mathcal{B(H)}$ and $\mu_1\mu_4=\mu_2\mu_3,$ then $$\|\mu_1PAP+\mu_2PAP^\bot+\mu_3P^\bot AP+\mu_4P^\bot AP^\bot\|\leq max\{|\mu_i|:\ i=1,2,3,4\}\|A\|,$$ where
$P^\bot=I-P.$

{\bf Proof.} We assume that $\mu_i\neq0$ for $i=1,2,3,4.$ Otherwise, the desired conclusion is clear. Since $\mu_1\mu_4=\mu_2\mu_3,$ it follows  that  $$(P+\frac{ \mu_3}{\mu_1}P^\perp)A(\mu_1P+\mu_2P^\perp)=\mu_1PAP+\mu_2PAP^\bot+\mu_3P^\bot AP+\mu_4P^\bot AP^\bot.$$
This implies  $$\begin{array}{rcl}\|\mu_1PAP+\mu_2PAP^\bot+\mu_3P^\bot AP+\mu_4P^\bot AP^\bot\|&\leq&\|P+\frac{ \mu_3}{\mu_1}P^\perp\|\|A\|\|(\mu_1P+\mu_2P^\perp)\|\\&=& max\{|\mu_i|:\ i=1,2,3,4\}\|A\|.\end{array} $$ $\Box$

{\bf Proposition 14.}  Let $\Psi:\mathcal{B(H)}\rightarrow \mathcal{B(H\otimes H)}$ be a linear map with the form  $$\Psi(X)=\lambda_1I \otimes X+\lambda_2X\otimes I+\lambda_3S(I\otimes X)+\lambda_4S(X\otimes I),$$ where $S\in\mathcal{B(H\otimes H)}$ is the swap operator and $\lambda_i$ are complex numbers for $i=1,2,3,4.$

 $(a)$ If $\lambda_1\lambda_2=\lambda_3\lambda_4,$ then $$\|\Psi\|_{cb}\leq max\{|\lambda_1+(-1)^j\lambda_2+(-1)^i(\lambda_3+(-1)^j\lambda_4)|:\ i,j=0,1\}.$$

$(b)$ If  $\lambda_1\lambda_2=\lambda_3\lambda_4$  and $\|\Psi(I)\|\geq max\{|\lambda_1-\lambda_2+(-1)^i(\lambda_3-\lambda_4)|:\ i=0,1\},$ then $\|\Psi\|_{cb}=\|\Psi\|=\|\Psi(I)\|.$

{\bf Proof.} $(a).$  Let $Q=\frac{I\otimes I+S}{2}$ and $Q^\perp=\frac{I\otimes I-S}{2}.$ Then \begin{equation} Q(I \otimes X)Q^\perp=\frac{1}{4}[I \otimes X-X\otimes I+S(I\otimes X)-S(X\otimes I)] \end{equation} and \begin{equation}   Q^\perp(I\otimes X)Q=\frac{1}{4}[I \otimes X-X\otimes I-S(I\otimes X)+S(X\otimes I)].\end{equation}
Combining equations (4.2)-(4.3) and (4.7)-(4.8), we get that \begin{equation}\begin{array}{rcl}\Psi(X)&=&\sum_{i=1}^{4}\lambda_iQ(I\otimes X)Q+(\lambda_1-\lambda_2+\lambda_3-\lambda_4)Q(I\otimes X)Q^\perp\\&&+(\lambda_1-\lambda_2-\lambda_3+\lambda_4)Q^\perp(I\otimes X)Q+\sum_{i=1}^{2}(\lambda_i-\lambda_{2+i})Q^\perp(I\otimes X)Q^\perp.\end{array} \end{equation}
For simplicity, we denote by $$\mu_1=\sum_{i=1}^{4}\lambda_i,\ \ \mu_2=\lambda_1-\lambda_2+\lambda_3-\lambda_4,\ \ \mu_3=\lambda_1-\lambda_2-\lambda_3+\lambda_4\ \hbox{ and }\ \mu_4=\sum_{i=1}^{2}(\lambda_i-\lambda_{2+i}).$$
Since $\lambda_1\lambda_2=\lambda_3\lambda_4,$ it follows that \begin{equation}\mu_1\mu_4=(\lambda_1+\lambda_2)^2-(\lambda_3+\lambda_4)^2=
(\lambda_1-\lambda_2)^2-(\lambda_3-\lambda_4)^2=\mu_2\mu_3.\end{equation}
For any $[A_{i,j}]_{n, n}\in\mathcal{M}_n({\mathcal{B(H)}}),$  we write $I\otimes A_{i,j}=B_{i,j}$ and $$\widetilde{Q}=diag(Q,Q,\cdots,Q)\in\mathcal{M}_n({\mathcal{B(H)}}).$$ Then $\|[A_{i,j}]_{n, n}\|=\|[B_{i,j}]_{n, n}\|$ and $\widetilde{Q}$ is an orthogonal projection with $$\widetilde{Q}^\bot=diag(Q^\bot,Q^\bot,\cdots,Q^\bot).$$ Thus we conclude from equations (4.9)-(4.10) and Lemma 13 that  $$\begin{array}{rcl}&&\|\Psi_n([A_{i,j}]_{n, n})\|\\&=&\|[\mu_1Q(I\otimes A_{i,j})Q+\mu_2Q(I\otimes A_{i,j})Q^\perp+\mu_3Q^\perp(I\otimes A_{i,j})Q+\mu_4Q^\perp(I\otimes A_{i,j})Q^\perp]_{n, n}\|\\&=&\|\mu_1 \widetilde{Q}[B_{i,j}]_{n, n}\widetilde{Q}+\mu_2 \widetilde{Q}[B_{i,j}]_{n, n}\widetilde{Q}^\perp+\mu_3 \widetilde{Q}^\perp[B_{i,j}]_{n, n}\widetilde{Q}+\mu_4 \widetilde{Q}^\perp[B_{i,j}]_{n, n}\widetilde{Q}^\perp\|\\&\leq&max\{|\mu_1|, |\mu_2|, |\mu_3|, |\mu_4|\}\|[B_{i,j}]_{n, n}\|\\&=&max\{|\mu_1|, |\mu_2|, |\mu_3|, |\mu_4|\}\|[A_{i,j}]_{n, n}\|,\end{array} $$ which implies that
$$\begin{array}{rcl}\|\Psi\|_{cb}=\sup_{n}\|\Psi_n\|&\leq& max\{|\mu_1|, |\mu_2|, |\mu_3|, |\mu_4|\}\\&=&max\{|\lambda_1+(-1)^j\lambda_2+(-1)^i(\lambda_3+(-1)^j\lambda_4)|:\ i,j=0,1\}.\end{array}$$

$(b).$ It is clear that $$\Psi(I)=(\lambda_1+\lambda_2)I\otimes I+(\lambda_3+\lambda_4)S,$$
so $$\|\Psi(I)\|=max\{|\lambda_1+\lambda_2+\lambda_3+\lambda_4|,\ |\lambda_1+\lambda_2-\lambda_3-\lambda_4|\}=max\{|\mu_1|,\ |\mu_4|\}.$$
If $\|\Psi(I)\|\geq max\{|\lambda_1-\lambda_2+(-1)^i(\lambda_3-\lambda_4)|:\ i=0,1\},$ then
$$max\{|\mu_1|,\ |\mu_4|\}\geq max\{|\mu_2|,\ |\mu_3|\}.$$ Thus it follows from $(a)$ that
$$\|\Psi(I)\|\leq\|\Psi\|\leq\|\Psi\|_{cb}\leq\|\Psi(I)\|,$$ and so $$\|\Psi\|_{cb}=\|\Psi\|=\|\Psi(I)\|$$ as desired. $\Box$

{\bf Remark 15.}  The assumption that $\Psi$ is a continuous linear map in the $W^*$-topology can be replaced by
the weak operator topology or the strong operator topology. Indeed, it is known that the set of all finite rank operators is also dense in $\mathcal{B(H)}$ with respect to the weak (strong) operator topology. Moreover,
for the sequence $X_n\in\mathcal{B(H)}$ and $X\in\mathcal{B(H)},$  $$\tau\hbox{-}\lim\limits_{n\rightarrow\infty}X_n=X\Longleftrightarrow \tau\hbox{-}\lim\limits_{n\rightarrow\infty}I\otimes X_n=I\otimes X\Longleftrightarrow \tau\hbox{-}\lim\limits_{n\rightarrow\infty}S(I\otimes X_n)=S(I\otimes X),$$ where $\tau$ is the $W^*$ or weak (strong) operator topology. $\Box$

\section{Proof of Theorem 4}

 We give the proof of $m=3.$ The proof technique of $m>3$ is analogous.
We conclude from  Schur-Weyl duality that if $\mathcal{H}$ is a finite dimensional Hilbert space and $T\in\mathcal{B(H\otimes H\otimes H)},$ then $$T(U\otimes U\otimes U)=(U\otimes U\otimes U)T$$ for all $U\in\mathcal{U(H)}$if and only if there exist complex numbers $\lambda_1,\lambda_2,\cdots,\lambda_6$  such that $T=\sum\limits_{i=1}^{6}\lambda_i\Gamma(s_{i}).$ Using the same technique as in the proof of Lemma 2, we  can show that the above conclusion holds for an infinite-dimensional Hilbert space ${\mathcal{H}}.$

 Let $$\Gamma(s_1)=I^{\otimes 3},\ \ \Gamma(s_2)=I\otimes S \ \hbox{ and } \
\Gamma(s_3)=S\otimes I,$$  where $S\in\mathcal{B(H\otimes H)}$ is the swap operator.
Moreover, we define $$\Gamma(s_4)(x_1\otimes x_2\otimes x_3)=x_3\otimes x_2\otimes x_1, \ \ \ \
\Gamma(s_5)(x_1\otimes x_2\otimes x_3)=x_3\otimes x_1\otimes x_2$$ and $$\Gamma(s_6)(x_1\otimes x_2\otimes x_3)=x_2\otimes x_3\otimes x_1$$ for all $x_1\otimes x_2\otimes x_3\in\mathcal{H}^{\otimes 3}.$

Suppose that $e\in \mathcal{H}$ is a fixed unit vector and
$$\mathcal{U}_e=\{ee^*+ V: \ V \hbox{ is a partial isometry with } ran(V)=ran(V^*)=ran(I-ee^*) \}.$$
For any $U\in \mathcal{U}_e,$ it is easy to see that $U$ is a unitary operator with $Uee^*=ee^*U,$ and so $$\Phi(ee^*)=\Phi(Uee^*U^*)=(U\otimes U\otimes U)\Phi(ee^*)(U^*\otimes U^*\otimes U^*).$$  It follows that $$[(ee^*+ V)\otimes (ee^*+ V)\otimes (ee^*+ V)]\Phi(ee^*)=\Phi(ee^*)[(ee^*+ V)\otimes (ee^*+ V)\otimes (ee^*+ V)].$$
Denote by $P=I-ee^*.$  Then using the same technique as in the proof of Theorem 2, we get that
\begin{equation}(ee^*\otimes ee^*\otimes ee^*)\Phi(ee^*)=\Phi(ee^*)(ee^*\otimes ee^*\otimes ee^*), \ \ (V\otimes V\otimes V)\Phi(ee^*)=\Phi(ee^*)(V\otimes V\otimes V), \end{equation}
\begin{equation}(ee^*\otimes V\otimes V+V\otimes V\otimes ee^*+V\otimes ee^*\otimes V) \Phi(ee^*)=\Phi(ee^*)(ee^*\otimes V\otimes V+V\otimes V\otimes ee^*+V\otimes ee^*\otimes V)\end{equation} and  \begin{equation}\begin{array}{rl}&(ee^*\otimes ee^*\otimes V+ee^*\otimes V\otimes ee^*+V\otimes ee^*\otimes ee^*) \Phi(ee^*)\\=&\Phi(ee^*)(ee^*\otimes ee^*\otimes V+ee^*\otimes V\otimes ee^*+V\otimes ee^*\otimes ee^*).\end{array}\end{equation} So
$ran(ee^*\otimes ee^*\otimes ee^*),$ $ran(P\otimes P\otimes P),$ $$ran(ee^*\otimes ee^*\otimes P)\oplus ran(ee^*\otimes P\otimes ee^*)\oplus ran(P\otimes ee^*\otimes ee^*)$$ and $$ran(ee^*\otimes P\otimes P)\oplus ran(P\otimes P\otimes ee^*)\oplus ran(P\otimes ee^*\otimes P)$$ are all reducing subspaces of $\Phi(ee^*).$
  Then equation (5.1) implies that there exist complex numbers $\mu$ and $\mu_{i1}$ (for $i=1,2,\cdots,6)$ such that \begin{equation}(ee^*\otimes ee^*\otimes ee^*)\Phi(ee^*)=(ee^*\otimes ee^*\otimes ee^*)\Phi(ee^*)(ee^*\otimes ee^*\otimes ee^*)=\mu(ee^*\otimes ee^*\otimes ee^*)\end{equation}
  and \begin{equation}(P\otimes P\otimes P)\Phi(ee^*)=(P\otimes P\otimes P)\Phi(ee^*)(P\otimes P\otimes P)=\sum\limits_{i=1}^{6}\mu_{i1}\Gamma(s_{i})(P\otimes P\otimes P).\end{equation}
 Moreover, we conclude from equation (5.2) that $$(ee^*\otimes V\otimes V)\Phi(ee^*)(P\otimes P\otimes ee^*)=(ee^*\otimes P\otimes P)\Phi(ee^*)(V\otimes V\otimes ee^*),$$ so $$\begin{array}{rcl}(V\otimes V\otimes ee^*)\Gamma(s_{6})\Phi(ee^*)(P\otimes P\otimes ee^*)&=&\Gamma(s_{6})(ee^*\otimes V\otimes V)\Phi(ee^*)(P\otimes P\otimes ee^*)\\&=&\Gamma(s_{6})(ee^*\otimes P\otimes P)\Phi(ee^*)(V\otimes V\otimes ee^*)\\&=&(P\otimes P\otimes ee^*)\Gamma(s_{6})\Phi(ee^*)(V\otimes V\otimes ee^*),\end{array}$$ which implies that there exist complex numbers $\nu_1$ and $\nu_2$ such that $$\begin{array}{rcl}(P\otimes P\otimes ee^*)\Gamma(s_{6})\Phi(ee^*)(P\otimes P\otimes ee^*)=\nu_1(P\otimes P\otimes ee^*)+\nu_2\Gamma(s_{3})(P\otimes P\otimes ee^*).\end{array}$$ Thus \begin{equation}\begin{array}{rcl}(ee^*\otimes P\otimes P)\Phi(ee^*)(P\otimes P\otimes ee^*)&=&\nu_1\Gamma(s_{6})(P\otimes P\otimes ee^*)+\nu_2\Gamma(s_{6})\Gamma(s_{3})(P\otimes P\otimes ee^*)\\&=&\nu_1\Gamma(s_{6})(P\otimes P\otimes ee^*)+\nu_2\Gamma(s_{2})(P\otimes P\otimes ee^*).\end{array}\end{equation}
 Analogously, we know that there exist complex numbers $\nu_i$ for $i=3,4,5,6$ such that  \begin{equation}\begin{array}{rcl}( P\otimes ee^*\otimes P)\Phi(ee^*)(P\otimes P\otimes ee^*)&=&\nu_3\Gamma(s_{5})(P\otimes P\otimes ee^*)+\nu_4\Gamma(s_{5})\Gamma(s_{3})(P\otimes P\otimes ee^*)\\&=&\nu_3\Gamma(s_{5})(P\otimes P\otimes ee^*)+\nu_4\Gamma(s_{4})(P\otimes P\otimes ee^*)\end{array}\end{equation}
 and \begin{equation}\begin{array}{rcl}( P\otimes P \otimes ee^*)\Phi(ee^*)(P\otimes P\otimes ee^*)=\nu_5 (P\otimes P\otimes ee^*)+\nu_6\Gamma(s_{3})(P\otimes P\otimes ee^*).\end{array}\end{equation}
 Then we conclude from equations (5.6)-(5.8) that there exist complex numbers  $\mu_{i4}$ for $i=1,2,\cdots,6$ such that $$(ee^*\otimes P\otimes P+ P \otimes ee^*\otimes P+P\otimes P\otimes ee^*)\Phi(ee^*)(P\otimes P\otimes ee^*)=\sum\limits_{i=1}^{6}\mu_{i4}\Gamma(s_{i})(P\otimes P\otimes ee^*).$$
 In a similar way, we get that there exist complex numbers  $\mu_{ij}$ (for $j=2,3$ and $i=1,2,\cdots,6)$ such that $$(ee^*\otimes P\otimes P+ P \otimes ee^*\otimes P+P\otimes P\otimes ee^*)\Phi(ee^*)( ee^*\otimes P\otimes P)=\sum\limits_{i=1}^{6}\mu_{i2}\Gamma(s_{i})(ee^*\otimes P\otimes P)$$ and $$(ee^*\otimes P\otimes P+ P \otimes ee^*\otimes P+P\otimes P\otimes ee^*)\Phi(ee^*)(P\otimes ee^*\otimes P)=\sum\limits_{i=1}^{6}\mu_{i3}\Gamma(s_{i})(P\otimes ee^*\otimes P).$$ Therefore,
  \begin{equation}\begin{array}{rl}&(ee^*\otimes P\otimes P+ P \otimes ee^*\otimes P+P\otimes P\otimes ee^*)\Phi(ee^*)\\=&\Phi(ee^*)(ee^*\otimes P\otimes P+ P \otimes ee^*\otimes P+P\otimes P\otimes ee^*)\\=&\sum\limits_{i=1}^{6}\mu_{i2}\Gamma(s_{i})(ee^*\otimes P\otimes P)+\sum\limits_{i=1}^{6}\mu_{i3}\Gamma(s_{i})(P\otimes ee^*\otimes P)+\sum\limits_{i=1}^{6}\mu_{i4}\Gamma(s_{i})(P\otimes P\otimes ee^*).\end{array}\end{equation}
 Using the same argument as in the proof of equation (5.9), we get from equation (5.3) that there exist complex numbers  $\mu_{ij}$ (for $j=5,6,7$ and  $i=1,2,\cdots,6)$ such that
  \begin{equation}\begin{array}{rl}&(ee^*\otimes ee^*\otimes P+ee^*\otimes P\otimes ee^*+ P \otimes ee^*\otimes ee^*)\Phi(ee^*)\\=&\Phi(ee^*)(ee^*\otimes ee^*\otimes P+ee^*\otimes P\otimes ee^*+ P \otimes ee^*\otimes ee^* )\\=&\sum\limits_{i=1}^{6}\mu_{i5}\Gamma(s_{i})(ee^*\otimes ee^*\otimes P)+\sum\limits_{i=1}^{6}\mu_{i6}\Gamma(s_{i})(ee^*\otimes P \otimes ee^*)+\sum\limits_{i=1}^{6}\mu_{i7}\Gamma(s_{i})(P\otimes ee^*\otimes ee^*).\end{array}\end{equation}
  Combining equations (5.4)-(5.5) and (5.9)-(5.10), we obtain that
    \begin{equation}\begin{array}{rl}\Phi(ee^*)&=\mu(ee^*\otimes ee^*\otimes ee^*)+\sum\limits_{i=1}^{6}\mu_{i1}\Gamma(s_{i})(P\otimes P\otimes P)+\sum\limits_{i=1}^{6}\mu_{i2}\Gamma(s_{i})(ee^*\otimes P\otimes P)\\&+\sum\limits_{i=1}^{6}\mu_{i3}\Gamma(s_{i})(P\otimes ee^*\otimes P)+\sum\limits_{i=1}^{6}\mu_{i4}\Gamma(s_{i})(P\otimes P\otimes ee^*)+\sum\limits_{i=1}^{6}\mu_{i5}\Gamma(s_{i})(ee^*\otimes ee^*\otimes P)\\&+\sum\limits_{i=1}^{6}\mu_{i6}\Gamma(s_{i})(ee^*\otimes P \otimes ee^*)+\sum\limits_{i=1}^{6}\mu_{i7}\Gamma(s_{i})(P\otimes ee^*\otimes ee^*)\\=&\lambda(ee^*\otimes ee^*\otimes ee^*)+\sum\limits_{i=1}^{6}\lambda_{i1}\Gamma(s_{i})(I\otimes I\otimes I)+\sum\limits_{i=1}^{6}\lambda_{i2}\Gamma(s_{i})(ee^*\otimes I\otimes I)\\&+\sum\limits_{i=1}^{6}\lambda_{i3}\Gamma(s_{i})(I\otimes ee^*\otimes I)+\sum\limits_{i=1}^{6}\lambda_{i4}\Gamma(s_{i})(I\otimes I\otimes ee^*)+\sum\limits_{i=1}^{6}\lambda_{i5}\Gamma(s_{i})(ee^*\otimes ee^*\otimes I)\\&+\sum\limits_{i=1}^{6}\lambda_{i6}\Gamma(s_{i})(ee^*\otimes I \otimes ee^*)+\sum\limits_{i=1}^{6}\lambda_{i7}\Gamma(s_{i})(I\otimes ee^*\otimes ee^*)
,\end{array}\end{equation}
 where $\lambda_{i1}=\mu_{i1},$ $\lambda_{ij}=\mu_{ij}-\mu_{i1}$ for $j=2,3,4,$  $$
  \lambda_{i5}=\mu_{i5}-\mu_{i2}-\mu_{i3}+\mu_{i1}, \ \ \lambda_{i6}=\mu_{i6}-\mu_{i2}-\mu_{i4}+\mu_{i1}, \  \ \  \lambda_{i7}=\mu_{i7}-\mu_{i3}-\mu_{i4}+\mu_{i1} $$ for $i=1,2,\cdots,6$  and $\lambda=\mu-\sum\limits_{i=1}^6\mu_{i1}-\sum\limits_{i=1}^6\sum\limits_{j=2}^{7}\lambda_{ij}.$

In a similar way to equation (2.33), we conclude that for any unit vector $f\in \mathcal{H},$  \begin{equation}\begin{array}{rl} \Phi(ff^*)&=\lambda(ff^*\otimes ff^*\otimes ff^*)+\sum\limits_{i=1}^{6}\lambda_{i1}\Gamma(s_{i})(I\otimes I\otimes I)+\sum\limits_{i=1}^{6}\lambda_{i2}\Gamma(s_{i})(ff^*\otimes I\otimes I)\\+&\sum\limits_{i=1}^{6}\lambda_{i3}\Gamma(s_{i})(I\otimes ff^*\otimes I)+\sum\limits_{i=1}^{6}\lambda_{i4}\Gamma(s_{i})(I\otimes I\otimes ff^*)+\sum\limits_{i=1}^{6}\lambda_{i5}\Gamma(s_{i})(ff^*\otimes ff^*\otimes I)\\+&\sum\limits_{i=1}^{6}\lambda_{i6}\Gamma(s_{i})(ff^*\otimes I \otimes ff^*)+\sum\limits_{i=1}^{6}\lambda_{i7}\Gamma(s_{i})(I\otimes ff^*\otimes ff^*).\end{array}\end{equation}

Let $ \widetilde{e}$ be a unit vector with $\widetilde{e}\bot e.$ Setting $g=\frac{1}{\sqrt{2}}e+\frac{1}{\sqrt{2}}\widetilde{e}$ and $h=\frac{1}{\sqrt{2}}e-\frac{1}{\sqrt{2}}\widetilde{e},$ we get that $$gg^*+hh^*=ee^*+\widetilde{e}\widetilde{e}^*,$$ so
\begin{equation}\Phi(gg^*)+\Phi(hh^*)=\Phi(ee^*)+\Phi(\widetilde{e}\widetilde{e}^*).\end{equation} Thus equations (5.12) and (5.13) imply that \begin{equation}\begin{array}{rl}
&\lambda(ee^*\otimes ee^*\otimes ee^*+\widetilde{e}\widetilde{e}^*\otimes\widetilde{e}\widetilde{e}^*\otimes\widetilde{e}\widetilde{e}^*) +\sum\limits_{i=1}^{6}\lambda_{i5}\Gamma(s_{i})(ee^*\otimes ee^*\otimes I+\widetilde{e}\widetilde{e}^*\otimes\widetilde{e}\widetilde{e}^*\otimes I)\\&+\sum\limits_{i=1}^{6}\lambda_{i6}\Gamma(s_{i})(ee^*\otimes I \otimes ee^*+\widetilde{e}\widetilde{e}^*\otimes I\otimes\widetilde{e}\widetilde{e}^*)+\sum\limits_{i=1}^{6}\lambda_{i7}\Gamma(s_{i})(I\otimes ee^*\otimes ee^*+I\otimes\widetilde{e}\widetilde{e}^*\otimes\widetilde{e}\widetilde{e}^*)\\=&\lambda(gg^*\otimes gg^*\otimes gg^*+hh^*\otimes hh^*\otimes hh^*)+\sum\limits_{i=1}^{6}\lambda_{i5}\Gamma(s_{i})(gg^*\otimes gg^*\otimes I+hh^*\otimes hh^*\otimes I)\\+&\sum\limits_{i=1}^{6}\lambda_{i6}\Gamma(s_{i})(gg^*\otimes I \otimes gg^*+hh^*\otimes I\otimes hh^*)+\sum\limits_{i=1}^{6}\lambda_{i7}\Gamma(s_{i})(I\otimes gg^*\otimes gg^*+I\otimes hh^*\otimes hh^*).\end{array}\end{equation} Since $dim\mathcal{H}\geq3,$ we can find a unit vector $\phi$ with $\phi\bot \bigvee\{e,\widetilde{e}\}.$ Then $\phi\bot g$ and $\phi\bot h,$ so $g\otimes h\otimes\phi\in\mathcal{H}\otimes\mathcal{H}\otimes\mathcal{H}.$ By equation (5.14), we have \begin{equation}\begin{array}{rcl}  \sum\limits_{i=1}^{6}\lambda_{i5}\Gamma(s_{i})(ee^*g\otimes ee^*h\otimes \phi+\widetilde{e}\widetilde{e}^*g\otimes\widetilde{e}\widetilde{e}^*h\otimes \phi)=0,\end{array}\end{equation}
  that is, \begin{equation}\begin{array}{rcl}   (\lambda_{15}+\lambda_{35})(e\otimes e\otimes \phi-\widetilde{e}\otimes\widetilde{e}\otimes \phi)&+&(\lambda_{25}+\lambda_{65})(e\otimes \phi\otimes e-\widetilde{e}\otimes \phi\otimes\widetilde{e})\\&+&(\lambda_{45}+\lambda_{55})(\phi\otimes e\otimes e- \phi\otimes \widetilde{e}\otimes\widetilde{e})=0.\end{array}\end{equation} Hence $$\lambda_{15}+\lambda_{35}=0,\ \ \ \lambda_{25}+\lambda_{65}=0 \ \hbox{ and }\ \lambda_{45}+\lambda_{55}=0.$$  For vectors $x,y,z\in \mathcal{H},$ we denote $\lambda_x=e^*x$ and $\lambda_y=e
^*y.$   Then \begin{equation}\begin{array}{rl}  &\sum\limits_{i=1}^{6}\lambda_{i5}\Gamma(s_{i})(ee^*\otimes ee^*\otimes I)(x\otimes y\otimes z)\\=& \lambda_x\lambda_y[(\lambda_{15}+\lambda_{35})(e\otimes e\otimes z)+(\lambda_{25}+\lambda_{65})(e\otimes z\otimes e)+(\lambda_{45}+\lambda_{55})(z\otimes e\otimes e)]=0.\end{array} \end{equation}
 Thus \begin{equation}\sum\limits_{i=1}^{6}\lambda_{i5}\Gamma(s_{i})(ee^*\otimes ee^*\otimes I)=0.\end{equation}
Analogously, we get that  \begin{equation}\sum\limits_{i=1}^{6}\lambda_{i6}\Gamma(s_{i})(ee^*\otimes I \otimes ee^*)=0 \ \ \hbox{ and }\ \  \sum\limits_{i=1}^{6}\lambda_{i7}\Gamma(s_{i})(I\otimes ee^*\otimes ee^*)=0.\end{equation}
 Then equation (5.14) yields that \begin{equation}\lambda(ee^*\otimes ee^*\otimes ee^*+\widetilde{e}\widetilde{e}^*\otimes\widetilde{e}\widetilde{e}^*\otimes\widetilde{e}\widetilde{e}^*)(g\otimes h\otimes g)=0,\end{equation} which implies \begin{equation}\lambda=0.\end{equation} Combining equations (5.11), (5.18)-(5.19) and (5.21), we have  \begin{equation}\begin{array}{rcl}\Phi(ee^*)&=& \sum\limits_{i=1}^{6}\lambda_{i1}\Gamma(s_{i})(I\otimes I\otimes I)+\sum\limits_{i=1}^{6}\lambda_{i2}\Gamma(s_{i})(ee^*\otimes I\otimes I) \\&+&\sum\limits_{i=1}^{6}\lambda_{i3}\Gamma(s_{i})(I\otimes ee^*\otimes I)+\sum\limits_{i=1}^{6}\lambda_{i4}\Gamma(s_{i})(I\otimes I\otimes ee^*).\end{array}\end{equation}
Thus equation (5.12) implies that for any unit vector $f\in \mathcal{H},$  \begin{equation}\begin{array}{rcl}\Phi(ff^*)&=&\sum\limits_{i=1}^{6}\lambda_{i1}\Gamma(s_{i})(I\otimes I\otimes I)+\sum\limits_{i=1}^{6}\lambda_{i2}\Gamma(s_{i})(ff^*\otimes I\otimes I)\\&+&\sum\limits_{i=1}^{6}\lambda_{i3}\Gamma(s_{i})(I\otimes ff^*\otimes I)+\sum\limits_{i=1}^{6}\lambda_{i4}\Gamma(s_{i})(I\otimes I\otimes ff^*).\end{array}\end{equation}
The remaining proof is the same to Theorem 1.

  Let $\{e_i\}_{i=1}^{m+1}$ be
an orthonormal basis of ${\mathcal{H}}.$ If there exist complex numbers $\mu_{ij}$  such that \begin{equation}\sum\limits_{i=1}^{m!}\sum\limits_{j=1}^{m+1}\lambda_{ij}\Gamma(s_{i})\Phi_j(X)=
\Phi(X)=\sum\limits_{i=1}^{m!}\sum\limits_{j=1}^{m+1}\mu_{ij}\Gamma(s_{i})\Phi_j(X)\end{equation} for all $X\in\mathcal{T(H)},$
then $$\sum\limits_{i=1}^{m!}\sum\limits_{j=2}^{m+1}(\lambda_{ij}-\mu_{ij})\Gamma(s_{i})\Phi_j(e_1e_2^*)=0,$$
  which implies $$\sum\limits_{i=1}^{m!}(\lambda_{i2}-\mu_{i2})\Gamma(s_{i})(e_1\otimes e_3\otimes\cdots\otimes e_{m+1})=\sum\limits_{i=1}^{m!}\sum\limits_{j=2}^{m+1}(\lambda_{ij}-\mu_{ij})\Gamma(s_{i})\Phi_j(e_1e_2^*)(e_2\otimes e_3\otimes\cdots\otimes e_{m+1})=0.$$ Hence $$\lambda_{i2}=\mu_{i2} \ \ \hbox{ for }\ i=1,2,\cdots,m!.$$ In a similar way, we get that $$\lambda_{ij}=\mu_{ij} \ \ \hbox{ for }j=3,4,\cdots,m+1 \hbox{ and } i=1,2,\cdots,m!.$$ Thus equation (5.24) implies that $$\sum\limits_{i=1}^{m!}\lambda_{i1}\Gamma(s_{i})\Phi_1(e_1e_1^*)
  =\sum\limits_{i=1}^{m!}\mu_{i1}\Gamma(s_{i})\Phi_1(e_1e_1^*),$$ and so $$\sum\limits_{i=1}^{m!}(\lambda_{i1}-\mu_{i1})\Gamma(s_{i})
  =0.$$ Then $$\lambda_{i1}=\mu_{i1} \ \ \hbox{ for }\ i=1,2,\cdots,m!.$$ $\Box$

%---------------------------------------------------------------------------------------%

\bibliographystyle{amsplain}

\begin{thebibliography}{99}
\bibitem {s1} I. Bardet, B. Collins, G. Sapra, Characterization of equivariant maps and application to entanglement detection. Ann. Henri Poincar\'{e}, 21 (2020), 3385-3406.
\bibitem {s1} B. V. R. Bhat,  Linear maps respecting unitary conjugation, Banach J. Math. Anal., 5 (2011), 1-5.
\bibitem {s1} B. Collins, H. Osaka, G. Sapra, On a family of linear maps from $M_n(C)$ to $M_{n^2}(C),$ Linear Algebra Appl., 555 (2018), 398-411.
\bibitem{s1}  M. D. Choi, Positive linear maps on C*-algebras, Can. Math. J., 24 (1972), 520-529.
\bibitem{s1} E. B. Davies, Quantum Theory of Open Systems, Academic Press, London, New York, San Francisco, (1976).
\bibitem {s1} J. Fullwood, General covariance for quantum states over time, Lett. Math. Phys., 114 (2024), 126.
\bibitem {s1} K. H. Han, S. H. Kye, E. St{\o}rmer, Infinite dimensional analogues of Choi matrices, J. Funct. Anal., 287 (2024), 110557.
\bibitem {s1}  J. Hou, A characterization of positive linear maps and criteria of entanglement for quantum states,
J. Phys. A, Math. Theor., 43 (2010), 385201.
\bibitem {s1}S. H. Kye, Positive maps in quantum information theory, Lecture Notes, Seoul National Univ,
2023, http://www.math.snu.ac.kr/~kye/book/qit.html.
\bibitem {s1} P. Kopszak, M. Mozrzymas and M. Studzi\'{n}ski, Positive maps from irreducibly
covariant operators, J. Phys. A: Math. Theor., 53 (2020), 395306.
 \bibitem {s1} I. Marvian, R. W. Spekkens, A generalization of Schur-Weyl duality
with applications in quantum estimation, Commun. Math. Phys., 331(2014), 431-475.
\bibitem {s1} Y. Li, H. K. Du, Interpolations of entanglement breaking channels and equivalent conditions for completely positive
maps, J. Funct. Anal., 268(2015), 3566-3599.
 \bibitem {s1} Y. Li, S. J. Wang, Completely positive and isometric maps on Schatten-class operators, Results Math., 80 (2025), Paper No. 41.
 \bibitem {s1} A. J. Parzygnat, J. Fullwood, F. Buscemi, G. Chiribella, Virtual quantum broadcasting, Phys. Rev. Lett., 132 (2024), 110203.
 \bibitem {s1} G. Popescu, Similarity and ergodic theory of positive linear maps, J. Reine
    Angew. Math., 561 (2003), 87-129.
 \bibitem {s1}   V. Paulsen, Completely bounded maps and operator algebras, Cambridge Studies in Advanced Mathematics,vol. 78, Cambridge University Press, Cambridge, 2002.
 \bibitem {s1} E. St{\o}rmer,  Positive Linear Maps of Operator Algebras, Springer-Verlag,
2013.
\bibitem {s1} E. St{\o}rmer, Separable states and positive maps, J. Funct. Anal., 254 (2008), 2303-2312.
 \bibitem {s1} K. Siudzi\'{n}ska, D. Chru\'{s}ci\'{n}ski, Quantum channels irreducibly covariant with respect to the finite group generated by the Weyl operators. J. Math. Phys. 59 (2018),  033508.
 \bibitem {s1} T. Takasaki and J. Tomiyama, On the geometry of positive maps in matrix
algebras, Math. Z., 184 (1983), 101-108.
\bibitem {s1} L. Zhang, Matrix integrals over unitary groups: An application of Schur-Weyl duality, arXiv: 1408.3782v5.

\end{thebibliography}

\end{document}